\documentclass[12pt]{article}
\usepackage{graphics}
\usepackage{amsmath}
\usepackage{amssymb}
\usepackage{graphicx}
\usepackage{makeidx}
\usepackage{mathrsfs}
\usepackage{amsfonts}
\usepackage[mathscr]{eucal}
\usepackage{amsmath}
\usepackage{mathrsfs}
\usepackage{epstopdf}
\usepackage{color}
\usepackage{geometry}
\def\ds{\displaystyle}
\def\RR{\vbox {\hbox to 8.9pt {I\hskip-2.1pt R\hfil}}}

\def\cen{\centerline}

\def\pni{\par\noindent}
\def\vsh{\smallskip}
\def\vs{\medskip}
\def\vvs{\bigskip}
\def\vsp{\vsh\pni} 

\def\ds{\displaystyle}

\newcommand{\bee}{\begin{equation}}
\newcommand{\ee}{\end{equation}}

\begin{document}

\cen{{\bf FRACALMO PRE-PRINT: \ http://www.fracalmo.org}}
\cen{{\bf Paper bublished in Mathematics (MDPI), Vol 8, 1273/1--27 (2021)}}
\cen{{\bf DOI:10.3390/math9111273}}

\vsh
\hrule
\vskip 0.50truecm
\font\title=cmbx12 scaled\magstep1
\font\bfs=cmbx12 scaled\magstep1
\font\little=cmr10
\begin{center}
{\title THE BATEMAN FUNCTIONS  
  REVISITED AFTER 90 YEARS:}
 \\ [0.25truecm]
 {\title A SURVEY OF OLD AND NEW RESULTS}
 \\ [0.45truecm]
{\bf Alexander Apelblat $^{1}$, Armando Consiglio $^{2}$, Francesco Mainardi $^{3}$}
\end{center}
$^{1}$  Department of Chemical Engineering, Ben Gurion University of the Negev, 84105 Beer Sheva, 84105, Israel.
\\
	$^{2}$  Institut f\"{u}r Theoretische Physik und Astrophysik and W\"{u}rzburg-Dresden Cluster of Excellence ct.qmat, Universit\"{a}t W\"{u}rzburg, D-97074 W\"{u}rzburg, Germany.
	\\
	$^{3}$ \ Dipartimento di  Fisica e Astronomia,  Universit\`{a} di Bologna, \& INFN, Via Irnerio 46, I-40126 Bologna, Italy.
\\
Corresponding author: francesco.mainardi@bo.infn.it
\vskip 0.5truecm
\centerline{\bf Revised Version: May  2021}

\begin{abstract}
\noindent 
The Bateman functions and the allied Havelock functions were introduced as solutions of some problems in hydrodynamics about ninety years ago, but after a period of one or two decades they were practically neglected. 
In handbooks, the Bateman function is only mentioned as a particular case of the confluent hypergeometric function. 
In order to revive our knowledge on these functions, their basic properties 
(recurrence functional and differential relations, series, integrals and the Laplace transforms) are presented. Some new results are also included. 
Special attention is directed to the Bateman and Havelock functions with integer orders, to  generalizations of these functions and to the Bateman-integral function
known in the literature
\end{abstract}
%

\vsp
{\it 2010 Mathematics Subject Classification (MSC)}:
  33C10, 33C15, 44A10.
\vsp
{\it Key Words and Phrases}: 
Bateman functions, Havelock functions, integral-Bateman functions, confluent hypergeometric functions.

\section{Introduction}
\noindent Harry Bateman (1882-1946) has been 
a renowned Anglo-American applied mathematician, who made outstanding contributions to mathematical physics, 
namely to aero- and fluid dynamics, to electro-magnetic and optical phenomena, to thermodynamics and geophysics and  to many other fields [1,2]. His main interests in mathematics were analytical solutions of partial differential and integral equations. His book published in 1932, \textit{Partial Differential Equations of Mathematical Physics} [3] is even today, a basic textbook on this subject. 
Born in Manchester and educated in Trinity College, Cambridge, with a continuation in 
Paris and Gottingen, Bateman emigrated to USA in 1910 and starting since 1917, during nearly three decades he has been Professor of Aeronautical Research and Mathematical Physics in the California Institute of Technology (Caltech). During these years he solved a number of various applied problems and simultaneously compiled from mathematical literature a vast amount of information associated with special functions and their properties. 

 From an enormous scientific legacy that Bateman left behind him, it is important to mention three items which are named after him. 
 The first is 
  the so-called  \textit{Bateman equation}
  which is applied in solutions of pbharmacokinetics problems (modeling of effective therapeutic management of drugs). 
 As  usual with Bateman, the origin of this equation came from an interaction with other scientists, and this  one  with Ernest Rutherford. It includes the solution of       a set of ordinary differential equations which describes the radioactive decay process. Mathematically, this process is similar to the behaviour of drugs in the human body and therefore is frequently used in pharmacokinetic models (see for example [4], and for prediction of the spread of COVID-19 look in [5]) listed in the fifties of the past century, and they constitute the so-called 
 {\it Bateman approach}.

 In mathematics, the Bateman name is mostly associated with 
the 
 five red books published in  the 
 fifties of the previous century, and they constitute
the  
  so-called  \textit{Bateman Manuscript Project}. 
  Three volumes are devoted to 
the  
  properties of special functions [1] and two volumes to tables of integral transforms [6]. 
  This enormous collection of functions, series and integrals, together with 
  the  
  description of their properties is based on the material compiled largely by Bateman, and prepared for publication by four editors A. Erd\'{e}lyi, R. Magnus, F. Oberhettinger and F.G. Tricomi. 
  Even today, these five books are indispensable for everybody, mathematicians, scientists and engineers who are involved in study and use of special functions and integral transforms. 
  They were essential as a precursor and model for later appearing in published or in modern on-line forms various compilations of mathematical reference data (for most important see for example [7-19]).

 In 1931 Bateman published a paper entitled: \textit{The $k$-function, a particular case of the confluent hypergeometric function}, where  he presented the definite trigonometric integral \eqref{GrindEQ__1_} and derived for it many properties [20]
\begin{equation} \label{GrindEQ__1_} 
\begin{array}{l} 
{k_{n} (x)=\frac{2}{\pi } \, \int _{0}^{\pi /2}\cos (x\tan \theta -n \theta )\, d\theta,  }
\quad  {n=0,1,2,3,...}
 \end{array} 
\end{equation} 
This integral represents the solution of 
 the
ordinary differential equation which appeared in Theodore von K\'{a}rm\'{a}n's theory of turbulent flows 
\begin{equation} \label{GrindEQ__2_} 
x \frac{d^{2} u(x)}{d x^{2} } =(x-n) u(x) \,.
\end{equation} 
Bateman named the integral in \eqref{GrindEQ__1_} as ${k}$-function in tribute for 
 the
outstanding contribution of von K\'{a}rm\'{a}n in the field of fluid dynamics. 
Nowadays, denoted in the mathematical literature by small or capital $k$, 
this function in a more general form, is called the Bateman function of argument 
$x$ and order (parameter) $\nu$. 
\begin{equation} \label{GrindEQ__3_} 
k_{\nu } (x)=\frac{2}{\pi } \, \int _{0}^{\pi /2}\cos (x\tan \theta -\nu  \theta )\, d\theta \,.  
\end{equation} 
The reason that Bateman used integer orders only, 
came from the fact that $k_{n}({x})$ functions can then be expressed by the Rodriguez type formulas and they are associated with the Laguerre polynomials. 
This also permitted to express sums of them in closed form and to link the Bateman functions with the confluent hypergeometric and Whittaker functions. 
In 1935, some new results were derived by Shastri [21], who showed that methods of operational calculus can be applied to this function.  

 Unfortunately, the Bateman functions found later rather limited attention in the mathematical literature. 
 Few only topics associated with them were considered and these mainly by Indian mathematicians [22-35]. They included the generalized Bateman functions, dual, triple and multi series equations of these functions, some integral equations and recurrence relations.
 \newpage
 
  It is worthwhile also to mention that in mathematical textbooks and tables, the Bateman function is not considered as a some kind of minor special function, but only indicated as a particular case of the confluent hypergeometric function. 
  Besides, no plots or tabulations of the Bateman functions are known in the literature. 

 One of the first attempts to enlarge a knowledge about properties of the Bateman functions, 
 has been evidently to introduce a new function, by replacing in the integrand of integral  \eqref{GrindEQ__1_} cosine function with sine function
\begin{equation} \label{GrindEQ__4_} 
 T_{n} (x)=\frac{2}{\pi } \, \int _{0}^{\pi /2}\sin (x\tan \theta -n \theta )\, d\theta, \quad n=0,1,2,3,... 
\end{equation} 
In 1950 H.M. Srivastava [25] and in 1966 K.N. Srivastava [29] suggested to denote this new function as $T_n(x)$, where the capital $T$ letter was adapted to honor Walter Tollmien who made pioneering works in the transition region between fully established laminar and turbulent flows. 
However,  an
unquestionably historical fact is that both trigonometric integrals as defined in \eqref{GrindEQ__1_} and \eqref{GrindEQ__4_}, were already, six year earlier in 1925, considered by Havelock who investigated some problems associated with surface waves [36].
 In the case of a circular cylinder immersed in a uniform flow, he needed to evaluated the following integrals which are written here in 
 their original notation for $k>0$
\begin{equation} \label{GrindEQ__5_} 
L_{r} =\int _{0}^{\pi /2}\cos (2  r\phi -k\tan \phi )\, d\phi  \,, \quad
M_{r} =\int _{0}^{\pi /2}\sin (2  r\phi -k\tan \phi )\, d\phi \,. 
\end{equation} 
Thus, in view of that $2r = x$  and ${k} = {n}$, these integrals differ from \eqref{GrindEQ__1_} and \eqref{GrindEQ__4_} only by the normalization 
factor $2/\pi$ and the minus sign in the second integral. 
What is even more important, Havelock was able to present the first six integrals in 
 a closed form. 
It is of interest also to mention that Bateman knew about the Havelock paper and of 
related integrals investigated by him. 
These integrals are included in the manuscript (later edited and published by 
Erd\'{e}lyi) which was found among his papers [37].
 Taking these facts into account, it is more fair and consistent to name the sine integral as the \textit{Havelock function} and to use similar as in \eqref{GrindEQ__3_} notation
\begin{equation} \label{GrindEQ__6_} 
h_{\nu } (x)=\frac{2}{\pi } \, \int _{0}^{\pi /2}\sin (x\tan \theta -\nu  \theta )\, d\theta  \,. 
\end{equation} 

 In the next step, further generalizations of the Bateman function were proposed by including powers of trigonometric functions in integrands for
 {$m,n=0,1,2,3,...$},
\begin{equation} \label{GrindEQ__7_} 
\begin{array}{l} 
{\ds k_{\nu }^{m} (x)}=
{\ds \frac{2}{\pi } \, \int _{0}^{\pi /2}(\cos \theta )^{m} \cos (x\tan \theta -\nu  \theta )\, d\theta  }\,,
 \\  \phantom{\rule{1pt}{15pt}}
 {\ds k_{\nu }^{m,n} (x)}=
 {\ds \frac{2}{\pi } \, \int _{0}^{\pi /2}(\sin \theta )^{m} (\cos \theta )^{n} \cos (x\tan \theta -\nu  \theta )\, d\theta  } \,. 
  \end{array} 
\end{equation} 
However, by reviewing the papers dealing with these so-called 
 generalized Bateman functions, Erd\'{e}lyi pointed out that the integrals in \eqref{GrindEQ__7_} are particular cases of confluent hypergeometric functions
  and the derived mathematical expressions are not new because they follow directly from manipulations with known properties of the Kummer confluent hypergeometric functions. 
\newpage

 Probably, the most paying attention from generalized Bateman functions is that which was proposed by Chaudhuri [38]. In an analogy with the integral Bessel functions, 
 he introduced the \textit{Bateman-integral function} 
\begin{equation} \label{GrindEQ__8_} 
ki_{n} (x)=-\int _{x}^{\infty }\frac{k_{2 n} (u)}{u}  \, du\quad ;\quad x>0, 
\end{equation} 
and discussed its properties. 

 As already mentioned above, in the last decades, the interest in the Bateman functions was very limited, and only investigations of Koepf and Schmersau [39-41] dealing with recurrence and other relations of  $F_n(x)$ functions, defined by 
\begin{equation} \label{GrindEQ__9_} 
\begin{array}{l} 
{\ds e^{-\, x (1\, +\, t)/(1\, -\, t)}}
 ={\ds \sum _{n\, =\, 0}^{\infty }t^{n}  F_{n} (x)},
  \\  \phantom{\rule{1pt}{15pt}}
  {\ds F_{n} (x)}
  =
  {\ds (-1)^{n} k_{2 n} (x)=(-1)^{n} \frac{2}{\pi } \, \int _{0}^{\pi /2}\cos (x\tan \theta -2 n \theta )\, d\theta \,. } 
  \end{array} 
\end{equation} 
should be mentioned.

 Considering that at the present time, the Bateman functions are unjustly neglected and nearly entirely  forgotten, we decided to prepare this survey in order to revive them and to promote them as independent functions. 
 It seems that the Bateman functions should be treated separately, less as particular cases of the confluent hypergeometric functions or the Whittaker functions. 
 Bearing in mind today that the literature on the subject is rather old and practically unknown, after Introduction, in the second section of this survey      we  collect  the most important properties of the Bateman functions with integer orders $k_{n}({x})$.
In the next section we present known results associated with the Havelock functions with integer orders $h_{n}({x})$.
  In the fourth section the generalized Bateman and Havelock functions
are discussed.   
 More general aspects related with the Bateman and Havelock functions having any order are considered in  the 
 fifth section. 
 In these sections some new results derived by us are also included. 
 The sixth section is dedicated to properties of the Bateman-integral functions.
  Concluding remarks are included in the last section. 
  
In Appendix A        we report  various finite and infinite integrals of functions associated with functions considered in this survey.
 Differential equations and trigonometric integrals associated with the Kummer confluent hypergeometric function are discussed in Appendix B.   
We refer the readers to  Appendix C
where they can find the  integral representations of known   special functions
recalled in the text because of their relations with the Bateman and Havelock
functions.
  
  It is expected that  all results presented here in analytical and in graphical form
  will stimulate a new research devoted to the Bateman and Havelock functions and these functions will find a desirable and proper place in the mathematical literature.

\vs

\section{The Bateman Functions with Integer  Orders}

\noindent The Bateman functions with integer order $n$ and with real argument 
$x$, are defined by
\begin{equation} \label{GrindEQ__10_} 
 k_{n} (x)=\frac{2}{\pi } \, \int _{0}^{\pi /2}\cos (x\tan \theta -n \theta )\, d\theta, \quad   {n=0,1,2,3,...}  
\end{equation} 
For this integral Bateman showed that [20] 
\begin{equation} \label{GrindEQ__11_}  
\begin{array}{l}
{\ds k_{n} (0)=\frac{2}{\pi  n} \, \sin \left(\frac{\pi  n}{2} \right)}, \quad
  {\ds k_{2 n} (0)=0}, 
\\   \phantom{\rule{1pt}{15pt}}
{\ds {\mathop{\lim }\limits_{x\, \to \, \infty }} k_{n} (x)={\mathop{\lim }\limits_{x\, \to \, \infty }} k'_{n} (x)=0},
   \end{array} 
\end{equation} 
and
\begin{equation} \label{GrindEQ__12_} 
\begin{array}{l}
 {\ds \left|k_{ n} (x)\right|\le 1}
  \\ \phantom{\rule{1pt}{15pt}}
   {\ds \left|k_{n} (x)\right|\le \left|\frac{n}{x} \right|\, ;\quad \left|k_{n} (x)\right|\le \left|\frac{n^{2} +2}{x^{2} } \right|\,;\quad n>2}, 
   \\  \phantom{\rule{1pt}{15pt}}
    {\ds \left|k_{2 n} (x)\right|\le \left|\frac{2 n}{x} \right|\ ;\quad x>1}, 
    \\  \phantom{\rule{1pt}{15pt}}
    {\ds \left|k'_{n} (x)\right|\le \left|\frac{n}{2 x} \right|}. 
    \end{array} 
\end{equation} 
In the case of even integers they are associated with the Havelock  integrals \eqref{GrindEQ__5_} and with $F_n(x)$ functions \eqref{GrindEQ__9_} 
in the following way [36,39-41] 
\begin{equation} \label{GrindEQ__13_} 
{k_{2n} (x)=\frac{2}{\pi } \, L_{n} (x)}, \quad
 {k_{2 n} (x)=(-1)^{n} F_{n} (x)}, \quad
 {h_{2n} (x)=-\frac{2}{\pi } \, M_{n} (x)}.
\end{equation} 
The first six Bateman functions were tabulated by Havelock[36]
for $x>0$,
\begin{equation} \label{GrindEQ__14_} 
\begin{array}{l} 
{\ds k_{0} (x)=e^{-\, x}}, 
\\  \phantom{\rule{1pt}{15pt}}
{\ds  k_{2} (x)=2 x e^{-\, x} },
 \\  \phantom{\rule{1pt}{15pt}}
  {\ds k_{4} (x)=2 x (x-1) e^{-\, x} },
   \\  \phantom{\rule{1pt}{15pt}}
   {\ds k_{6} (x)=\frac{2}{3}  x (2 x^{2} -6 x+3) e^{-\, x} }, 
   \\  \phantom{\rule{1pt}{15pt}}
    {\ds k_{8} (x)=\frac{2}{3}  x (x^{3} -6 x^{2} +9 x-3) e^{-\, x} },
     \\  \phantom{\rule{1pt}{15pt}}
     {\ds k_{10} (x)=\frac{2}{15}  x (2 x^{4} -20 x^{3} +60 x^{2} -60 x+15) e^{-\, x} },
      \\  \phantom{\rule{1pt}{15pt}}
 {\ds k_{12} (x)=
 \frac{2}{45}  x (2 x^{5} -30 x^{4} +150 x^{3} -300 x^{2} +225 x-45) e^{-\, x} } .
    \end{array} 
\end{equation} 
In    the
general case these polynomials can be derived from the Rodriguez type formula
\begin{equation} \label{GrindEQ__15_} 
k_{2n} (x)=\frac{(-1)^{n}  x e^{x} }{n!} \, 
\frac{d^{n} }{d x^{n} } \, \left[x^{n -1} e^{- 2 x} \right] ,
\end{equation} 
which is similar to that of the generalized Laguerre polynomials $L_{n}^{(\alpha )} (x)$. 
\begin{equation} \label{GrindEQ__16_} 
L_{n}^{(\alpha )} (x)=\frac{ x^{-\alpha }  e^{x} }{n!} \, \frac{d^{n} }{d x^{n} } \, 
\left[x^{n + \alpha } e^{- x} \right]. 
\end{equation} 
Bateman showed that for his functions with even integer orders we have [20]
\begin{equation} \label{GrindEQ__17_} 
k_{2n} (x)=(-1)^{n}   e^{-\, x} \left[L_{n} (2 x)-L_{n\, -\, 1} (2 x)\right], 
\end{equation} 
where $L_{k}({z})$ are the Laguerre polynomials..

  It is more difficult to express the Bateman functions with odd orders in terms of other known functions. For $n = 1$, Bateman introduced a new integration variable 
  $ {t} = \tan \theta$   and obtained [20] 
\begin{equation} \label{GrindEQ__18_} 
\begin{array}{l} 
{\ds k_{1} (x)=\frac{2}{\pi } \, \int _{0}^{\pi /2}\cos (x\tan \theta -\theta )\, d\theta = } 
\\  \phantom{\rule{1pt}{15pt}}
 {\ds \frac{2}{\pi } \, \int _{0}^{\pi /2}\cos (x\tan \theta )\cos \theta \, d\theta +\frac{2}{\pi } \, \int _{0}^{\pi /2}\sin (x\tan \theta )\sin \theta \, d\theta =  }
  \\  \phantom{\rule{1pt}{15pt}}
  {\ds \frac{2}{\pi } \, \int _{0}^{\infty }\frac{\cos (x t)}{(1+t^{2} )^{3/2} } \, dt
  + \frac{2}{\pi } \, \int _{0}^{\infty }\frac{t \sin (x t)}{(1+t^{2} )^{3/2} } \, dt =} 
  \\  \phantom{\rule{1pt}{15pt}}
   {\ds \frac{2}{\pi } \, \int _{0}^{\infty }\frac{\cos (x t)}{(1+t^{2} )^{3/2} } \, dt
   - \frac{2 x}{\pi } \, \int _{0}^{\infty }\frac{\cos (x t)}{(1+t^{2} )^{1/2} } \, dt. } 
   \end{array} 
\end{equation} 
The last two integrals are the integral representations  of the modified Bessel functions of the second kind of the first and zero orders [7]
\begin{equation} \label{GrindEQ__19_} 
\begin{array}{l} 
{\ds k_{1} (x)=\frac{2 x}{\pi } \, \left[K_{1} (x)-K_{0} (x)\right]\, ;\quad x>0}, 
\\  \phantom{\rule{1pt}{15pt}}
{\ds  k_{1} (x)=-\frac{2 x}{\pi } \, \left[K_{1} (-\, x)+K_{0} (-\, x)\right]\,;\quad x<0}. \end{array} 
\end{equation} 
The Bateman functions with other even and odd integer orders can also be derived by applying the recurrence relations which are in the form of difference equations and differential-difference equations 
\begin{equation} \label{GrindEQ__20_} 
\begin{array}{l} 
{\ds (2 x-2 n)\, k_{2 n} (x)=(n-1)\, k_{2 n\, -\, 2} (x)+(n+1)\, k_{2 n\, +\, 2} (x)}
 \\  \phantom{\rule{1pt}{15pt}}
  {\ds 4 x  k'_{ n} (x)=(n-2)\, k_{ n\, -\, 2} (x)-(n+2)\, k_{ n\, +\, 2} (x)} 
  \\  \phantom{\rule{1pt}{15pt}}
  {\ds k'_{ n} (x)+k'_{ n\, +\, 2} (x)=k_{n} (x)-k_{ n\, +\, 2} (x)}
   \\ \phantom{\rule{1pt}{15pt}}
    {\ds x k''_{ n} (x)=(x-n)\, k_{n} (x).}
   \end{array} 
\end{equation} 
For example, using the second equation in \eqref{GrindEQ__20_} for ${n} = 1$,
 we have
\begin{equation} \label{GrindEQ__21_} 
\begin{array}{l} 
{\ds k_{3} (x)=-\frac{1}{3} \left[4 x\frac{d k_{1} (x)}{dx} +k_{-\, 1} (x) \right]} 
\\  \phantom{\rule{1pt}{15pt}}
 {\ds \frac{d k_{1} (x)}{dx} =\frac{2}{\pi } \left[K_{1} (x)-K_{0} (x)\right]+\frac{2 x}{\pi } \, \left[\frac{d K_{1} (x)}{dx} -\frac{d K_{0} (x)}{dx} \right]}
  \\  \phantom{\rule{1pt}{15pt}}
  {\ds \frac{d K_{1} (x)}{dx} =\frac{K_{2} (x)+K_{0} (x)}{2} } 
  \\ \phantom{\rule{1pt}{15pt}}
  {\ds \frac{d K_{0} (x)}{dx} =-\, K_{1} (x)}
   \end{array} 
\end{equation} 
and $k_{-1}({x})$  can be expressed by using integrals from \eqref{GrindEQ__18_}
\begin{equation} \label{GrindEQ__22_} 
\begin{array}{l} 
{\ds k_{-\, 1} (x)=\frac{2}{\pi } \, \int _{0}^{\pi /2}\cos (x\tan \theta +\theta )\, d\theta = }
 \\  \phantom{\rule{1pt}{15pt}}
  {\ds \frac{2}{\pi } \, \int _{0}^{\pi /2}\cos (x\tan \theta )\cos \theta \, d\theta -\frac{2}{\pi } \, \int _{0}^{\pi /2}\sin (x\tan \theta )\sin \theta \, d\theta.   }
   \end{array} 
\end{equation} 

 It is also possible to obtain the Bateman functions with odd orders in 
 a different new procedure, for example $k_{3}(x)$ 
\begin{equation} \label{GrindEQ__23_} 
\begin{array}{l}
 {\ds k_{3} (x)=\frac{2}{\pi } \, \int _{0}^{\pi /2}\cos (x\tan \theta -3 \theta )\, d\theta = }
  \\  \phantom{\rule{1pt}{15pt}}
   {\ds \frac{2}{\pi } \, \int _{0}^{\pi /2}\cos (x\tan \theta )\cos (3 \theta )\, d\theta 
   +\frac{2}{\pi } \, \int _{0}^{\pi /2}\sin (x\tan \theta )\sin (3 \theta )\, d\theta,   }
    \end{array} 
\end{equation} 
but with ${t=\tan \theta}$ 
\begin{equation} \label{GrindEQ__24_} 
\begin{array}{l} 
{\ds \sin (3 \theta )=3\sin \theta -4 (\sin \theta )^{3} 
=
\sin \theta \, \frac{3-(\tan \theta )^{2} }{1+(\tan \theta )^{2} } 
=
\frac{t (3-t^{2} )}{(1+t^{2} )^{3/2} } },
 \\  \phantom{\rule{1pt}{15pt}}
 {\ds \cos (3 \theta )=-3\cos \theta +4 (\cos \theta )^{3}
  =
  \cos \theta \, \frac{1-3(\tan \theta )^{2} }{1+(\tan \theta )^{2} } 
  =\frac{(1-3 t^{2} )}{(1+t^{2} )^{3/2} }, } 
   \end{array} 
\end{equation} 
and therefore \eqref{GrindEQ__23_} becomes
\begin{equation} \label{GrindEQ__25_} 
k_{3} (x)=\frac{2}{\pi } \, \int _{0}^{\infty }\frac{(1-3 t^{2} )\, \cos (xt)}{(1+t^{2} )^{5/2} } \, dt+\frac{2}{\pi } \, \int _{0}^{\infty }\frac{t (3-t^{2} )\, \sin (xt)}{(1+t^{2} )^{5/2} } \, dt.   
\end{equation} 
However, this type of integrals can be evaluated by differentiating the modified Bessel functions of the second kind [14]   
\begin{equation} \label{GrindEQ__26_} 
\begin{array}{l} 
{\ds \int _{0}^{\infty }\frac{t^{2 n\, +\, 1} \, \sin (xt)}{(1+t^{2} )^{\alpha } } \, dt
 =(-1)^{n\, +\, 1} \frac{2^{1/2\, -\, \alpha } \sqrt{\pi } }{\Gamma (\alpha )} \, \frac{\partial ^{2 n +\, 1} }{\partial  x^{2 n +\, 1} } \, \left[x^{\alpha \, -\, 1/2}  K_{\alpha \, -\, 1/2} (x)\right]}, 
\quad {\ds \alpha >n+1/2 },
 \\  \phantom{\rule{1pt}{15pt}}
  {\ds  \int _{0}^{\infty }\frac{t^{2 n} \, \sin (xt)}{(1+t^{2} )^{\alpha } } \, dt =(-1)^{n} \frac{2^{1/2\, -\, \alpha } \sqrt{\pi } }{\Gamma (\alpha )} \, \frac{\partial ^{2 n +\, 1} }{\partial  x^{2 n +\, 1} } \, \left[x^{\alpha \, -\, 1/2}  K_{\alpha \, -\, 1/2} (x)\right]},
  \quad {\ds \alpha >n}.
    \end{array} 
\end{equation} 
Using known expressions for
$\sin(\alpha+ 2\theta)$
 and 
 $\cos(\alpha+ 2\theta)$ functions with
  $\alpha  = 2{n} + 1$, and taking into 
  account that [7] with ${t=\tan \theta }$
\begin{equation} \label{GrindEQ__27_} 
\begin{array}{l}
 {\ds \sin (2 \theta )=\frac{2 \tan \theta }{1+(\tan \theta )^{2} } =\frac{2 t }{(1+t^{2} )} },
  \\  \phantom{\rule{1pt}{15pt}}
  {\ds \cos (2 \theta )=\, \frac{1-(\tan \theta )^{2} }{1+(\tan \theta )^{2} } =\frac{(1- t^{2} )}{(1+t^{2} )}, }
   \end{array} 
\end{equation} 
the above described procedure can be extended to the Bateman functions with higher odd  orders. 
Integrals of the type presented in \eqref{GrindEQ__26_} can be also used when derivatives with respect to the argument are considered
with ${m=0,1,2,3,...}$ 
\begin{equation} \label{GrindEQ__28_} 
\begin{array}{l} 
{\ds \frac{\partial ^{2 m} \, k_{n} (x)}{\partial  x^{2 m} } 
=(-1)^{m} \frac{2}{\pi } \,
 \int _{0}^{\pi /2}(\tan \theta )^{2 m} \cos (x\tan \theta -n \theta )\, d\theta,  } 
 \\  \phantom{\rule{1pt}{15pt}}
  {\ds \frac{\partial ^{2 m\, +\, 1} \, k_{n} (x)}{\partial  x^{2 m\, +\, 1} }
   =(-1)^{m} \frac{2}{\pi } \, \int _{0}^{\pi /2}(\tan \theta )^{2 m\, +\, 1}
    \sin (x\tan \theta -n \theta )\, d\theta.  } 
     \end{array} 
\end{equation} 

 In order to illustrate the behaviour of  the Bateman functions as a function of argument and order, they were numerically evaluated using the MATLAB program and they are presented in Figure 1 for positive integer orders and in Figure 2 for negative integer order. 
 As can be observed by comparing both figures, the curves are shifted with the symmetry predicted by Bateman [20]
\begin{equation} \label{GrindEQ__29_} 
k_{-\, n} (x)=k_{n} (-\, x)\,. 
\end{equation}

\begin{figure}[h!]
	\centering
	\includegraphics[width=12cm]{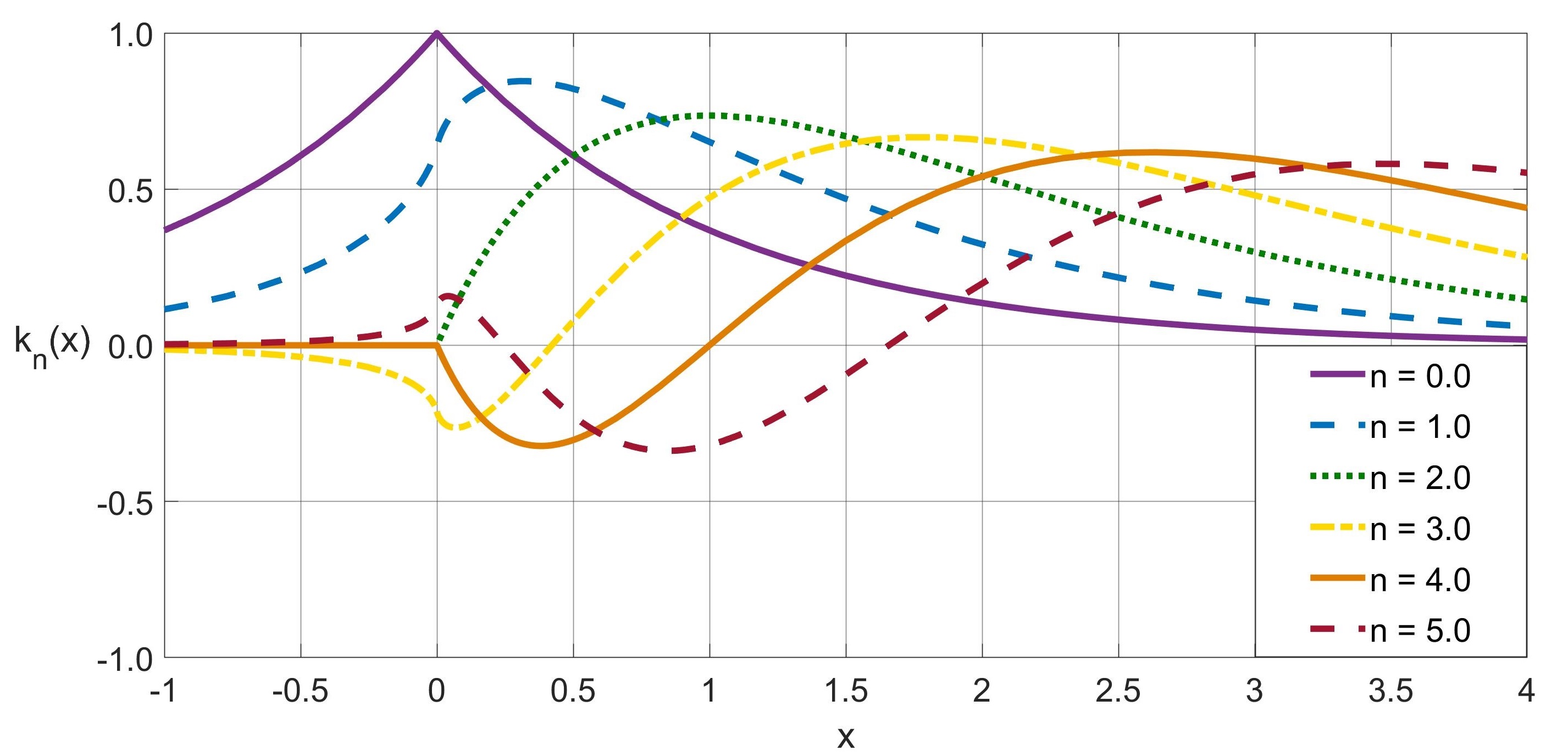}
	\caption{Bateman functions with  positive integer orders as a function of argument \textit{x}.} 
\end{figure}

\begin{figure}[h!]
	\centering
	\includegraphics[width=12cm]{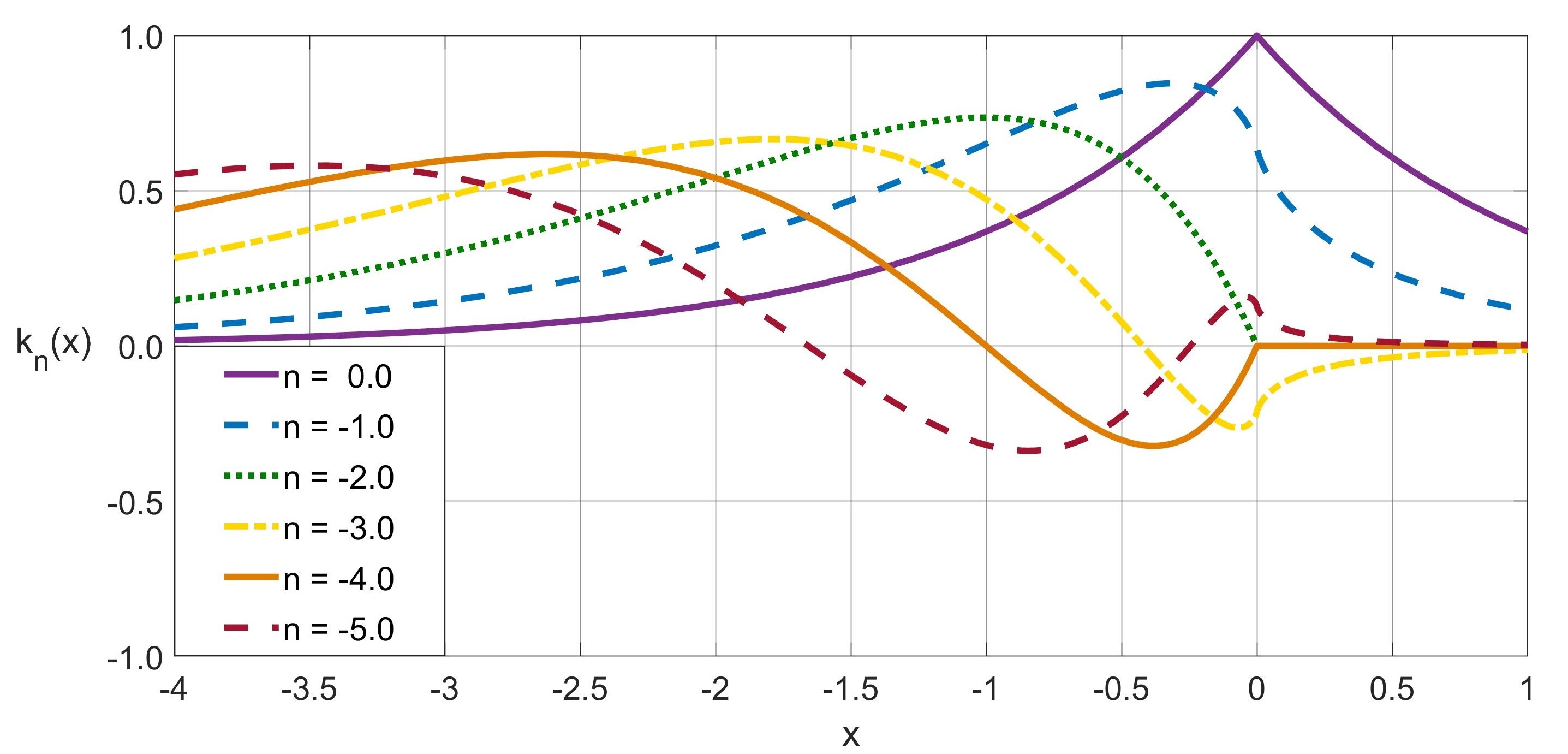}
	\caption{Bateman functions with  negative integer orders as a function of argument \textit{x}.} 
\end{figure}
\noindent
Considering similarity with the generalized Laguerre polynomials, Bateman was able to show 
     the
existence of the following expansions associated with his functions with even orders [20] 
\begin{equation} \label{GrindEQ__30_} 
\begin{array}{l} 
{\ds \sum _{n\, =\, 0}^{\infty }(-1)^{n}   t^{n}  k_{2 n} (x)
=
(1-t)^{\alpha \, +\, 1} e^{-\, x} \sum _{n\, =\, 0}^{\infty }   t^{n}  L_{n}^{(\alpha )} (2 x)},
 \\  \phantom{\rule{1pt}{15pt}}
  {\ds \sum _{n\,  =\, 0}^{\infty }\frac{t^{n} }{2^{n}  n!} k_{2 n\, +\, 2} (x)
   =
   2 e^{-\, (x\, +t/2)} \sqrt{\frac{x}{t} } \, I_{1} (2\sqrt{x t} )},
   \\  \phantom{\rule{1pt}{15pt}} 
   {\ds \sum _{n = 0}^{\infty }(-1)^{n}   k_{4 n +2} (x)=\sin x},
\quad   
  {\ds \sum _{n\, =\, 0}^{\infty }(-1)^{n}   k_{4 n} (x)=\cos x}.,
  \end{array} 
\end{equation} 
where $I_1$ denotes the modified Bessel function of order 1, see (C.8) and [7].
Shabde [22] demonstrated that 
\begin{equation} \label{GrindEQ__31_} 
\begin{array}{l} 
{\ds \sum _{n\, =\, 0}^{\infty }(n+1)  t^{n}  k_{2 n\, +\, 2} (x) =\frac{2 x e^{-\, x\, +[2 x t/(1+t)]} }{(1+t)^{2} } },
 \\  \phantom{\rule{1pt}{15pt}} 
 {\ds \sum _{n\, =\, 0}^{\infty }(-1)^{n} (2 n+1)  t^{2 n}  k_{2 n\, +\, 2} (x) =} 
 \\  \phantom{\rule{1pt}{15pt}}
 {\ds \frac{2 x e^{-\, x\, +2 x t^{2} /(1+t^{2} )} }{(1+t^{2} )^{2} } \left[(1-t^{2} ) \cos \left(\frac{2 x t}{1+t^{2} } \right)+2 t \sin \left(\frac{2 x t}{1+t^{2} } \right)\right]}, 
 \\   \phantom{\rule{1pt}{15pt}}
 {\ds \sum _{n\, =\, 0}^{\infty }(-1)^{n} (2 n+2)  t^{2 n\, +\, 1}  k_{4 n\, +\, 4} (x) =} 
 \\  \phantom{\rule{1pt}{15pt}}
 {\ds \frac{2 x e^{-\, x\, +2 x t^{2} /(1+t^{2} )} }{(1+t^{2} )^{2} } \left[(1-t^{2} ) \sin \left(\frac{2 x t}{1+t^{2} } \right)-2 t \cos \left(\frac{2 x t}{1+t^{2} } \right)\right]}, \end{array} 
\end{equation} 
and
\begin{equation} \label{GrindEQ__32_} 
\begin{array}{l}
 {\ds \sum _{n\, =\, 0}^{\infty }\frac{(-1)^{n}  t^{n} }{  n!} k_{2 n\, +\, 2} (x) =\sqrt{\frac{2 x}{t} } \, e^{-\, (x\, +t)} \, J_{1} (2^{3/2} \sqrt{x t} )},
  \\  \phantom{\rule{1pt}{15pt}}
   {\ds \sum _{n\, =\, 0}^{\infty }\frac{(-1)^{ n}  t^{2 n} }{(2 n)!} k_{2 n\, +\, 2} (x) =\sqrt{\frac{2 x}{t} } \, \left[-\sin t\, ber'(2^{3/2} \sqrt{x t} )+\cos t\, bei'(2^{3/2} \sqrt{x t} )\right]},
    \\  \phantom{\rule{1pt}{15pt}}
    {\ds \sum _{n\, =\, 0}^{\infty }\frac{(-1)^{ n\, +\, 1}  t^{2 n + 1 } }{(2 n+1)!} k_{4 n\, +\, 4} (x) =\sqrt{\frac{2 x}{t} } \, \left[\cos t\, ber'(2^{3/2} \sqrt{x t} )+\sin t\, bei'(2^{3/2} \sqrt{x t} )\right]},
     \end{array} 
\end{equation} 
where $ber'(z)$ and $bei'(z)$  are 
the
derivatives of the Kelvin functions.

\noindent Additional sums of series expansions were reported by Shastri [24]
\begin{equation} \label{GrindEQ__33_} 
\begin{array}{l} 
{\ds \sum _{n\, =\, 0}^{\infty }(-1)^{n}  t^{2 n\, +\, 1}  k_{4 n\, +\, 2} (x) =\, e^{ x (t^{2} \, -\, 1)/(1+t^{2} )} \sin \left(\frac{2 x t}{1+t^{2} } \right)\, ;\quad \left|t\right|<1}, 
\\  \phantom{\rule{1pt}{15pt}}
{\ds \sum _{n\, =\, 0}^{\infty }(-1)^{n}  t^{2 n}  k_{4 n} (x) =\, e^{ x (t^{2} \, -\, 1)/(1+t^{2} )} \cos \left(\frac{2 x t}{1+t^{2} } \right)\, ;\quad \left|t\right|<1}, 
\\   \phantom{\rule{1pt}{15pt}}
{\ds \sum _{n\, =\, 0}^{\infty }(-1)^{n}  k_{4 n\, +\, 2} (x) =\sin x},
 \quad  {\ds \sum _{n\, =\, 0}^{\infty }(-1)^{n}  k_{4 n} (x) =\cos x},
  \end{array} 
\end{equation} 
and
\begin{equation} \label{GrindEQ__34_} 
\begin{array}{l} 
{\ds \sum _{n\, =\, 0}^{\infty }k_{2 n} (x) \sin (2 n \theta )=\sin (x\tan \theta )},
 \\ \phantom{\rule{1pt}{15pt}}
  {\ds \sum _{n\, =\, 0}^{\infty }k_{2 n} (x) \sin (2 n \theta )=\sin (x\tan \theta )}, 
  \\  \phantom{\rule{1pt}{15pt}}
  {\ds \sum _{n\, =\, 0}^{\infty }k_{2 n} (x) =1}. 
  \end{array} 
\end{equation} 

 The orthogonal relations were established by Bateman [20]
\begin{equation} \label{GrindEQ__35_} 
\begin{array}{l} 
{\ds \int _{0}^{\infty }[k_{2 n} (x)]^{2} \,  dx=\left\{\begin{array}{c} {1\quad ;\quad n>0} 
\\  \phantom{\rule{1pt}{15pt}}
{1/2\quad ;\quad n=0} \end{array}\right. }
 \\  \phantom{\rule{1pt}{25pt}}
  {\ds \int _{0}^{\infty }k_{2 n} (x) k_{2 n\, +\, 2 k} (x) \, dx=\left\{\begin{array}{c} {0\quad ;\quad k>0} 
  \\ \phantom{\rule{1pt}{15pt}}
  {1/2\quad ;\quad k=0} \end{array}\right. }
   \\  \phantom{\rule{1pt}{25pt}}
    {\ds \int _{0}^{\infty }\frac{k_{n} (x) k_{2 k} (x)}{x}  \, dx=\frac{4 \sin \left[\frac{\pi }{2}  (2 k-n)\right]}{\pi  n (2 k-n)} \, ;\quad k>0}, 
    \end{array} 
\end{equation} 
and over the entire integration interval 
\begin{equation} \label{GrindEQ__36_} 
\begin{array}{l} 
{\ds \int _{-\, \infty }^{+\infty }k_{2 k} (x) k_{2 m} (x) \, dx=\frac{\sin [\pi  (m-k)]}{\pi  (k-m+1)\, (k-m)\, (k-m-1)} },
 \\  \phantom{\rule{1pt}{15pt}}
 {\ds PV\int _{-\, \infty }^{\infty }k_{2 k\, +\, 1} (x) k_{2 m\, +\, 1} (x) \, \frac{dx}{x} =
 \left\{
 \begin{array}{c} {0\, ;\quad k\ne m,}
 \\ {\frac{2}{\pi  (2k+1)} \, ;\quad k=m.} 
 \end{array}
 \right. }
  \end{array} 
\end{equation} 

 In the literature there is a number of infinite integrals where the Bateman functions appear in integrands or in final results of integration. 
 These integrals are collected in Appendix A, here only the Laplace transforms of the Bateman functions  are presented [6,9]:
\begin{equation} \label{GrindEQ__37_} 
\begin{array}{l} 
{\ds \int _{0}^{\infty }e^{-\, s t} k_{0} (t)\,  dt=\frac{1}{s+1} \, ;
\quad Re(s+1)>0\,;\quad n=0,1,2,...}
 \\  \phantom{\rule{1pt}{15pt}}
  {\ds \int _{0}^{\infty }e^{-\, s t} k_{2 n\, +\, 2} ( t)\,  dt=\frac{2  (1-s)^{n} }{(s+1)^{n\, +\, 2} } }
   \\  \phantom{\rule{1pt}{15pt}}
    {\ds \int _{0}^{\infty }e^{-\, s t} k_{2 \nu } ( t)\,  dt=\frac{\sin (\pi  \nu )}{2 \pi  \nu  (1-\nu ) } \, _{2} F_{1} (1,2;2-\nu ;\frac{1-s}{2} )\, ;\quad Res>0} 
    \end{array} 
\end{equation} 
and
\begin{equation} \label{GrindEQ__38_} 
\begin{array}{l} 
{\ds \int _{0}^{\infty }e^{-\, s t} e^{-\, t^{2} } k_{2 n} (t^{2} )\,  dt=\frac{(-1)^{n\, -\, 1} s^{n\, -\, 3/2} \, e^{s^{2} /16} }{2^{3 n/2\, +\, 1/4} } \, W_{-\, n/2\, -1/4,-\, n/2\, -1/4} \left(\frac{s^{2} }{8} \right)}
 \\  \phantom{\rule{1pt}{15pt}}
  {\ds \int _{0}^{\infty }e^{-\, s t} k_{2 m\, +\, 2} (\frac{t}{2} )\, k_{2 n\, +\, 2} (\frac{t}{2} )\,  \frac{dt}{t} =\frac{(-1)^{m\, +\, n} }{(s+1)^{m\, +\, n\, +\, 2} } \, _{2} F_{1} \left(-m,-n;2;\frac{1}{s^{2} } \right)}
   \\  \phantom{\rule{1pt}{15pt}}
    {\ds Res>-1} 
    \\  \phantom{\rule{1pt}{15pt}}
     {\ds \int _{0}^{\infty }e^{-\, s t} \frac{e^{(\alpha \, +\, \beta ) t/2} }{\alpha  \beta } k_{2 m\, +\, 2} (\frac{\alpha  t}{2} )\, k_{2 n\, +\, 2} (\frac{\beta  t}{2} )\,  \frac{dt}{t} =}
      \\  \phantom{\rule{1pt}{15pt}}
       {\ds \frac{(-1)^{m\, +\, n} (m+n+1)!\, (s-\alpha )^{m} \, (s-\beta )^{n} }{(m+1)!\, (n+1)!(s+1)^{m\, +\, n\, +\, 2} } \, _{2} F_{1} \left(-m,-n;-m-n-1;\frac{s (s-\alpha -\beta )}{(s-\alpha )\, (s-\beta )} \right)} 
       \\  \phantom{\rule{1pt}{15pt}}
       {\ds m,n=0,1,2,...\, ;\quad Res>0}
        \end{array} 
\end{equation} 
where $W_{\kappa, \mu}(z)$  is the Whittaker function. 
Formulas in (32 and (33) are accessible in a much more general forms by applying the basic properties of the Laplace transformation
\begin{equation} \label{GrindEQ__39_} 
\begin{array}{l} 
{\ds L\left\{f(t)\right\}=\int _{0}^{\infty }e^{-\, s t} f(t)\,  dt=F(s)\, ;\quad a>0} 
\\  \phantom{\rule{1pt}{15pt}}
 {\ds L\left\{f(a t)\right\}=\frac{1}{a} \, F\left(\frac{s}{a} \right)}
  \\  \phantom{\rule{1pt}{15pt}}
  {\ds L\left\{e^{\pm \, a t} f( t)\right\}=F\left(s\mp a\right)}
   \\  \phantom{\rule{1pt}{15pt}}
    {\ds L\left\{t^{n}  f( t)\right\}=(-1)^{n} \frac{d^{n}  F(s)}{d s^{n} } }
     \end{array} 
\end{equation} 
For example in the simple case of  the function $k_2(t)$ 
 we have from \eqref{GrindEQ__39_}
\begin{equation} \label{GrindEQ__40_} 
\begin{array}{l} 
{\ds L\left\{k_{2} (t)\right\}=\frac{2}{(s+1)^{2} } }
 \\  \phantom{\rule{1pt}{15pt}}
  {\ds L\left\{k_{2} (a t)\right\}=\frac{2 a}{(s+a)^{2} } } 
  \\  \phantom{\rule{1pt}{15pt}}
  {\ds L\left\{e^{\pm \, a t} k_{2} (a t)\right\}=\frac{2 }{(s\mp a+1)^{2} } } 
  \\  \phantom{\rule{1pt}{15pt}}
  {\ds L\left\{t k_{2} (t)\right\}=\frac{4}{(s+1)^{3} } }.
   \end{array} 
\end{equation} 
The initial and final values of the Bateman functions 
with even integer orders (see Figure 1) as presented in \eqref{GrindEQ__11_}, 
can also be derived from 
  the
rules of the operational calculus
\begin{equation} \label{GrindEQ__41_} 
\begin{array}{l}
 {\ds k_{0} (t\, \to \, +0)={\mathop{\lim }\limits_{s\, \to \, \infty }} [s F(s)]={\mathop{\lim }\limits_{s\, \to \, \infty }} \left[\frac{s}{s+1} \right]=1} 
 \\  \phantom{\rule{1pt}{15pt}}
  {\ds k_{0} (t\, \to \, \infty )={\mathop{\lim }\limits_{s\, \to \, 0}} [s F(s)]={\mathop{\lim }\limits_{s\, \to \, 0}} \left[\frac{s}{s+1} \right]=0} 
  \\  \phantom{\rule{1pt}{15pt}}
   {\ds k_{2 n\, +\, 2} (t\, \to \, +0)={\mathop{\lim }\limits_{s\, \to \, \infty }} [s F(s)]={\mathop{\lim }\limits_{s\, \to \, \infty }} \left[\frac{2  s (1-s)^{n} }{(s+1)^{n\, +\, 2} } \right]=0} 
   \\  \phantom{\rule{1pt}{15pt}}
   {\ds k_{2 n\, +\, 2} (t\, \to \, \infty )={\mathop{\lim }\limits_{s\, \to \, 0}} [s F(s)]={\mathop{\lim }\limits_{s\, \to \, 0}} \left[\frac{2  s (1-s)^{n} }{(s+1)^{n\, +\, 2} } \right]=0}.
   \end{array} 
\end{equation} 
Since the Bateman function is a particular case of the Whittaker function 
\begin{equation} \label{GrindEQ__42_} 
k_{2 \nu } \left(\frac{t}{2} \right)=\frac{1}{\Gamma (\nu +1)} \, W_{\nu ,1/2} (t) 
\end{equation} 
it is possible to enlarge a number of the Laplace transforms using transforms of the Whittaker functions $W_{1/2, 1/2}(t)$  and
$W_{\nu, 1/2}(t)$
\begin{equation} \label{GrindEQ__43_} 
\begin{array}{l}
 {\ds L\left\{t^{1/2}  e^{1/2 t} k_{1} \left(\frac{2}{t} \right)\right\}=\frac{\sqrt{\pi } }{s} \, \left[H_{1} (2 \sqrt{s} )-Y_{1} (2 \sqrt{s} )\right]} 
 \\  \phantom{\rule{1pt}{15pt}}
 {\ds L\left\{t e^{1/2 t} k_{1} \left(\frac{2}{t} \right)\right\}=\frac{1}{2 s} \, H_{1}^{(1)} (\sqrt{s} )\, H_{1}^{(2)} (\sqrt{s} )}
  \\  \phantom{\rule{1pt}{15pt}}
  {\ds L\left\{\frac{1}{t}   e^{-\, 1/2 t} k_{1} \left(\frac{2}{t} \right)\right\}=\frac{2^{5/2}  \sqrt{s} }{\pi } \, K_{0} (\sqrt{s} )\, K_{1} (\sqrt{s} )} 
  \\  \phantom{\rule{1pt}{15pt}}
  {\ds L\left\{\frac{1}{t^{2} }   e^{-\, 1/2 t} k_{1} \left(\frac{2}{t} \right)\right\}=\frac{4}{\pi  s} \, \left[K_{1} (\sqrt{s} )\right]^{2} }
   \end{array} 
\end{equation} 
and
\begin{equation} \label{GrindEQ__44_} 
\begin{array}{l}
 {\ds L\left\{t^{\alpha \, -\, 1}   k_{2 \nu } \left(\frac{t}{2} \right)\right\}=} 
 \\  \phantom{\rule{1pt}{15pt}}
 {\ds \frac{\Gamma (\alpha )\, \Gamma (\alpha +1)}{\Gamma (\nu +1) \Gamma (\alpha -\nu +1)} \left(\frac{2}{2 s+1} \right)^{\alpha \, +\, 1} \, _{2} F_{1} (\alpha +1,-\, \nu ;\alpha -\nu +1;\frac{2 s-1}{2 s+1} )}
  \\  \phantom{\rule{1pt}{15pt}}
  {\ds Res>-\frac{1}{2} }
   \\  \phantom{\rule{1pt}{15pt}}
    {\ds L\left\{t^{\nu }   e^{\, 1/2 t} k_{2 \nu } \left(\frac{2}{t} \right)\right\}=\frac{2^{1\, -2 \nu } }{\Gamma (\nu +1)\, s^{\nu \, +\, 1/2} } S_{2 \nu ,1} (2 \sqrt{s} )\, ;
    \quad Re(\nu \pm \frac{1}{2} )>-\frac{1}{2} }
     \\  \phantom{\rule{1pt}{15pt}}
     {\ds L\left\{\frac{1}{t^{\nu } }   e^{-\, 1/2 t} k_{2 \nu } \left(\frac{2}{t} \right)\right\}=\frac{2 s^{\nu \, -\, 1/2} }{\Gamma (\nu +1)} K_{1} (2 \sqrt{s} )\, ;\quad Res>0,} 
     \end{array} 
\end{equation} 
where
$H_\mu(t)$,
$Y_\mu(t)$,
$H_\mu^{(1)}(t)$, $H_\mu^{(2)}(t)$
and $S_\mu(t)$
 are  the 
  Struve, Bessel, Hankel and Lommel functions, respectively.

\vs

\section{The Havelock Functions with Integer  Orders}

\noindent As pointed out above, Havelock in solving the surface wave problem [36] 
encountered the following trigonometric integrals with even integer values of order (parameter) $n$
\begin{equation} \label{GrindEQ__45_} 
h_{n} (x)=\frac{2}{\pi } \, \int _{0}^{\pi /2}\sin (x\tan \theta -n \theta )\, d\theta.   
\end{equation} 
These functions with positive and negative values of order were calculated numerically by using the MATLAB program and they are plotted in Figures 3 and 4. 
Comparing both figures, it is evident that the curves are shifted according to
\begin{equation} \label{GrindEQ__46_} 
h_{-\, n} (x)=-\, h_{n} (-\, x). 
\end{equation} 
Havelock was able to present the first six integrals in terms of polynomials and the logarithmic integrals [36]
\begin{equation} \label{GrindEQ__47_} 
\begin{array}{l} 
{\ds h_{0} (x)=\frac{1}{2} \left[e^{x}  li(e^{-\, x} )-e^{-\, x}  li(e^{x} )\right] }
 \\  \phantom{\rule{1pt}{15pt}}
 {\ds h_{2} (x)=x e^{-\, x}  li(e^{x} )-1} 
 \\  \phantom{\rule{1pt}{15pt}}
 {\ds h_{4} (x)=x (x-1) e^{-\, x}  li(e^{x} )-x} 
 \\ \phantom{\rule{1pt}{15pt}}
 {\ds h_{6} (x)=\frac{x (2 x^{2} -6 x+3) e^{-\, x}  li(e^{x} )-(2 x^{2} -4 x+1)}{3}  }
  \end{array} 
\end{equation} 
and

\begin{figure}[h!]
	\centering
	\includegraphics[width=12cm]{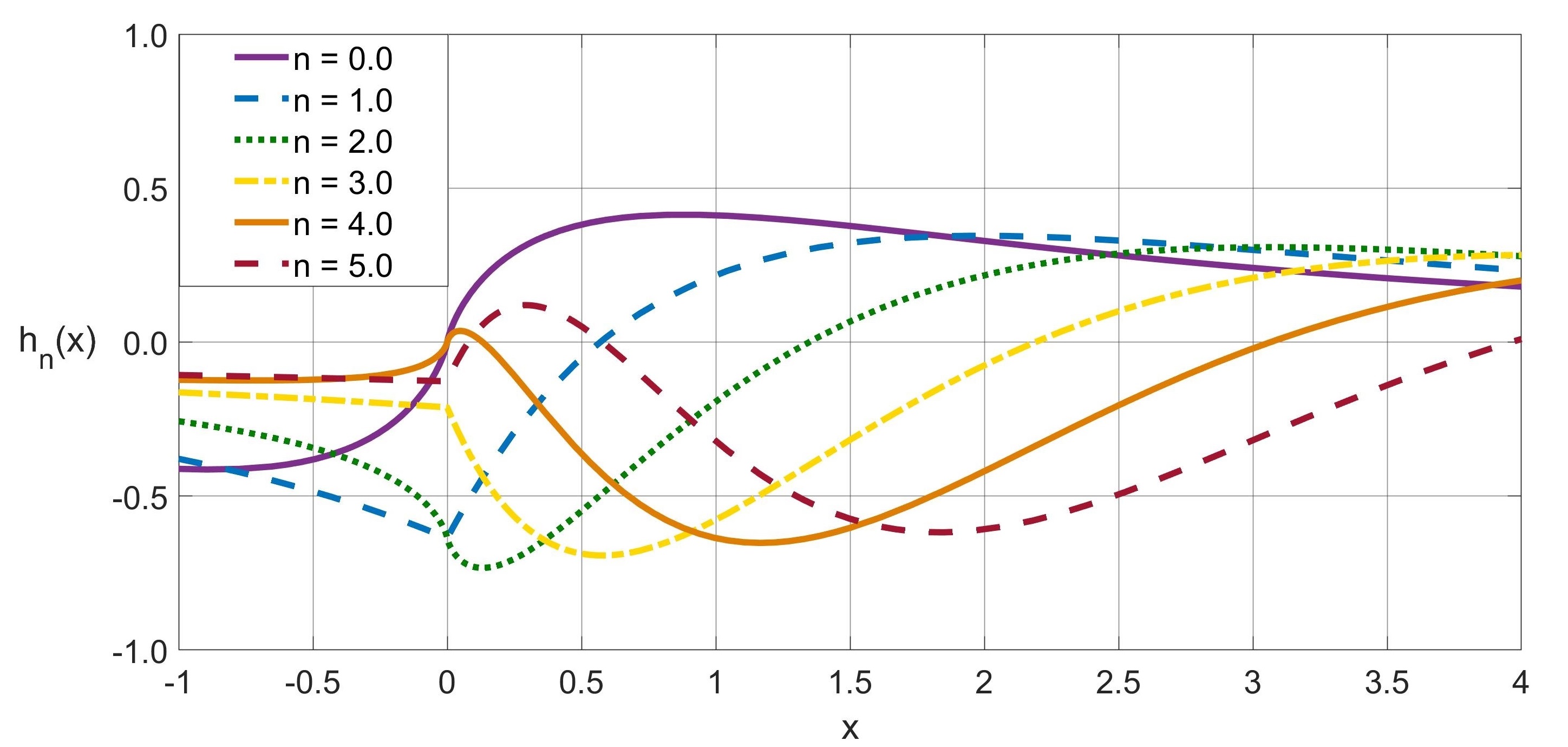}
	\caption{Havelock  functions with  positive integer orders as a function of argument \textit{x}.} 
\end{figure}

\begin{figure}[h!]
	\centering
	\includegraphics[width=12cm]{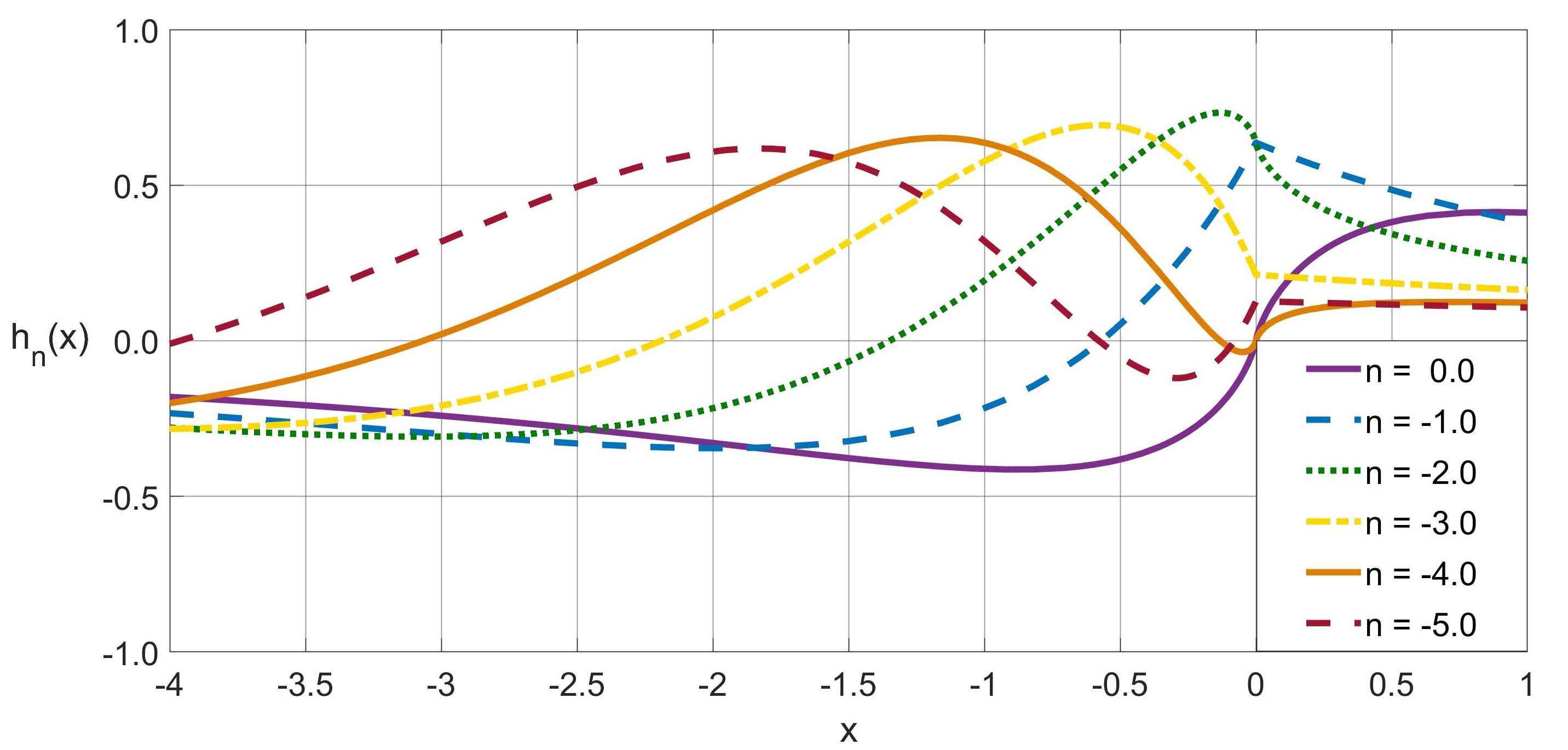}
	\caption{Havelock  functions with begative integer orders as a function of argument \textit{x}.} 
\end{figure}

\begin{equation} \label{GrindEQ__48_} 
\begin{array}{l} 
{\ds h_{8} (x)=\frac{x (x^{3} -6 x^{2} +9 x-3) e^{-\, x}  li(e^{x} )-x  (x^{2} -5 x+5)}{3}  }
 \\  \phantom{\rule{1pt}{15pt}}
  {\ds h_{10} (x)=\frac{x (2 x^{4} -20 x^{3} +60 x^{2} -60 x+15) e^{-\, x}  li(e^{x} )}{15}  -}
   \\ \phantom{\rule{1pt}{15pt}} 
  {\ds \frac{(2 x^{4} -18 x^{3} +44 x^{2} -28 x+3)}{15} }
   \\  \phantom{\rule{1pt}{15pt}}
    {\ds h_{12} (x)=\frac{x (2 x^{5} -30 x^{4} +150 x^{3} -300x^{2} +225 x-45) e^{-\, x}  li(e^{x} )}{45}  -} 
    \\  \phantom{\rule{1pt}{15pt}} 
    {\ds \frac{x  (2 x^{4} -28 x^{3} +124 x^{2} -198x+93)}{45} } 
    \end{array} 
\end{equation} 
where
\begin{equation} \label{GrindEQ__49_} 
{li(z)=\int _{0}^{z}\frac{d t}{\ln t}  =\gamma +\ln z+\sum _{n\, =\, 1}^{\infty }\frac{z^{n} }{n!n}  }, \quad   {z=e^{x} }. 
\end{equation} 
In the same way as in the Bateman paper from 1931, 
 the
properties of the Havelock functions with integer orders were studied by Srivastava in 1950 [25]. He found that
\begin{equation} \label{GrindEQ__50_} 
\begin{array}{l} 
{\ds \left|h_{n} (x)\right|\le 1} 
\\  \phantom{\rule{1pt}{15pt}}
 {\ds h_{n} (0)=\frac{2}{\pi  n} \left[\cos \left(\frac{\pi  n}{2} \right)-1\right]}
  \\  \phantom{\rule{1pt}{15pt}}
   {\ds h_{2 n} (0)=\left[\frac{1-(-1)^{n} }{\pi  n} \right]} \\ {h_{4 n} (0)=0} 
   \\  \phantom{\rule{1pt}{15pt}}
    {\ds {\mathop{\lim }\limits_{x\, \to \, \infty }} h_{n} (x)={\mathop{\lim }\limits_{x\, \to \, \infty }} h_{n} {}^{{'} } (x)=0},
     \end{array} 
\end{equation} 
and
\begin{equation} \label{GrindEQ__51_} 
\begin{array}{l}
 {\ds h_{0} (x)=\frac{2}{\pi } \, \int _{0}^{\pi /2}\sin (x\tan \theta )\, d\theta  =\frac{2}{\pi } \, \int _{0}^{\infty }\frac{\sin (x t)}{1+t^{2} } \, dt }
  \\  \phantom{\rule{1pt}{15pt}}
   {\ds h_{1} (x)=\frac{2}{\pi } \, \int _{0}^{\pi /2}\sin (x\tan \theta -\theta )\, d\theta  =} 
   \\  \phantom{\rule{1pt}{15pt}}
    {\ds \frac{2}{\pi } \, \int _{0}^{\pi /2}\left[\sin (x\tan \theta )\, \cos \theta -\cos (x\tan \theta )\, \sin \theta \right]\, d\theta  =}
     \\  \phantom{\rule{1pt}{15pt}}
      {\ds \frac{2}{\pi } \, \int _{0}^{\infty }\frac{[\sin (x t)-t  \cos (x t)]}{(1+t^{2} )^{3/2} } \, dt }. 
      \end{array} 
\end{equation} 
These integrals are of the type presented in \eqref{GrindEQ__26_}. 
In 1950 Srivastava [25] showed that the infinite integral in \eqref{GrindEQ__51_} can be expressed in terms of
 the
 modified Bessel function of the first kind of zero order and the Struve function of zero order and their derivatives. 
\\
 The Havelock  functions satisfy the following recurrence and differential relations [25,37] 
\begin{equation} \label{GrindEQ__52_} 
\begin{array}{l} 
{\ds (2 n-4 x)\, h_{n} (x)+(n-2)\, h_{ n\, -\, 2} (x)+(n+2)\, h_{ n\, +\, 2} (x)=-\frac{8}{\pi } }
 \\  \phantom{\rule{1pt}{15pt}}
 {\ds 4 x h'_{n} (x)=(n-2)\, h_{ n\, -\, 2} (x)-(n+2)\, h_{ n\, +\, 2} (x)} 
 \\  \phantom{\rule{1pt}{15pt}}
  {\ds h'_{ n\, -\, 1} (x)+h'_{ n\, +\, 1} (x)=h_{n\, -\, 1} (x)-h_{n\, +\, 1} (x)}
   \\  \phantom{\rule{1pt}{15pt}}
    {\ds x h''_{ n} (x)=(x-n)\, h_{n} (x)-\frac{2}{\pi } }. 
    \end{array} 
\end{equation} 
The Laplace transform of the function $h_{0}(x)$ can be  obtained in the following way
\begin{equation} \label{GrindEQ__53_} 
\begin{array}{l} 
{\ds L\left\{h_{0} (x)\right\}=\frac{2}{\pi } \, \int _{0}^{\infty }e^{-\, s x}  \, \left(\int _{0}^{\pi /2}\sin (x\tan \theta )\, d\theta  \right)\, dx=} 
\\   \phantom{\rule{1pt}{15pt}}
{\ds \frac{2}{\pi } \, \int _{0}^{\pi /2}\left(\, \int _{0}^{\infty }e^{-\, s x} \sin (x\tan \theta )\, dx \right) \, d\theta =
\frac{2}{\pi } \, \int _{0}^{\pi /2}\frac{\tan \theta }{s^{2} +(\tan \theta )^{2} } \, d\theta = }
 \\  \phantom{\rule{1pt}{15pt}}
  {\ds \frac{2}{\pi } \, \int _{0}^{\infty }\frac{t}{(s^{2} +t^{2} ) (1+t^{2} )}  \, dt=\frac{2\, \ln (s)}{\pi  (s^{2} -1)} }.
   \end{array} 
\end{equation} 
For the function $h_{1}(x)$   we have
\begin{equation} \label{GrindEQ__54_} 
\begin{array}{l} 
{\ds L\left\{h_{1} (x)\right\}=\frac{2}{\pi } \, \int _{0}^{\pi /2}\left(\, \int _{0}^{\infty }e^{-\, s x} \left[\sin (x\tan \theta )\, \cos \theta -\cos (x\tan \theta )\, \sin \theta \right]dx \right) \, d\theta =} 
\\  \phantom{\rule{1pt}{15pt}}
 {\ds =\frac{2 }{\pi } \, \int _{0}^{\pi /2}\frac{[\tan \theta \cos \theta -s\, \sin \theta ]}{s^{2} +(\tan \theta )^{2} } \, d\theta = \frac{2 (1-s)}{\pi } \, \int _{0}^{\infty }\frac{t}{(s^{2} +t^{2} ) (1+t^{2} )^{3/2} }  \, dt=}
  \\  \phantom{\rule{1pt}{15pt}}
  {\ds \frac{2 }{\pi  (s+1)} \, \left[\frac{\sec ^{-1} (s)}{\sqrt{s^{2} -1} } -1\right]}.
   \end{array} 
\end{equation} 
The Laplace transforms of the functions $h_0(x)$ and $h_1(x)$   were also derived by Srivastava [25] in 1950 , but in the final expressions, the factor $2/\pi$ is missing.

 The Havelock function $h_2(x)$  is expressed by
\begin{equation} \label{GrindEQ__55_} 
\begin{array}{l}
 {\ds h_{2} (x)=\frac{2}{\pi } \, \int _{0}^{\pi /2}\sin (x\tan \theta -2 \theta )\, d\theta  =}
  \\  \phantom{\rule{1pt}{15pt}}
   {\ds \frac{2}{\pi } \, \int _{0}^{\pi /2}\left[\sin (x\tan \theta )\, \cos (2 \theta )-\cos (x\tan \theta )\, \sin (2 \theta )\right]\, d\theta  =}
    \\  \phantom{\rule{1pt}{15pt}}
     {\ds \frac{2}{\pi } \, \int _{0}^{\pi /2}\frac{\sin (x\tan \theta )\, [1-(\tan \theta )^{2} ]-2 \tan \theta \cos (x\tan \theta )\, }{1+(\tan \theta )^{2} } d\theta  }
      \\  \phantom{\rule{1pt}{15pt}}
       {\ds \frac{2}{\pi } \, \int _{0}^{\infty }\frac{(1-t^{2} )\sin (x t)-2 t  \cos (x t)}{(1+t^{2} )^{2} } \, dt }
        \end{array} 
\end{equation} 
and its Laplace transform is therefore
\begin{equation} \label{GrindEQ__56_} 
\begin{array}{l} 
{\ds L\left\{h_{2} (x)\right\}=\frac{2}{\pi } \, \int _{0}^{\infty } \left[\int _{0}^{\infty }e^{-\, s x}  \frac{(1-t^{2} )\sin (x t)-2 t  \cos (x t)}{(1+t^{2} )^{2} } \, dx\right]\, dt =}
 \\  \phantom{\rule{1pt}{15pt}}
  {\ds -\frac{2  [s+1+\ln (s)]}{\pi  (s+1)^{2} } }
   \end{array} 
\end{equation} 
where the infinite integrals in \eqref{GrindEQ__53_}, \eqref{GrindEQ__54_} and \eqref{GrindEQ__56_} were verified using 
 the
MATHEMATICA program. 
The derived Laplace transforms allow us  to obtain the initial and final values
 of the Havelock  functions,  for example for the function $h_0(x)$
  we have
\begin{equation} \label{GrindEQ__57_} 
\begin{array}{l}
 {\ds h_{0} (x\, \to \, +0)={\mathop{\lim }\limits_{s\, \to \, \infty }} [s F(s)]={\mathop{\lim }\limits_{s\, \to \, \infty }} \left[\frac{2 s\ln (s)}{s^{2} -1} \right]=0} 
 \\  \phantom{\rule{1pt}{15pt}}
  {\ds h_{0} (x\, \to \, \infty )={\mathop{\lim }\limits_{s\, \to \, 0}} [s F(s)]={\mathop{\lim }\limits_{s\, \to \, 0}} \left[\frac{2 s\ln (s)}{s^{2} -1} \right]=0}
   \end{array} 
\end{equation} 
as it is observed in Figure 3.

 There is a number of recurrence and differential expressions that include both the Bateman and the Havelock  functions. 
 They were reported by Srivastava [25] and three of them are presented here
\begin{equation} \label{GrindEQ__58_} 
\begin{array}{l} 
{\ds (n-2)\, \left[k_{n} (x) h_{ n\, -\, 2} (x)-k_{ n\, -\, 2} (x) h_{n} (x)\right]+}
 \\   \phantom{\rule{1pt}{15pt}}
 {\ds (n+2)\, \left[k_{n} (x) h_{ n\, +\, 2} (x)-k_{ n\, +\, 2} (x) h_{n} (x)\right]=-\frac{8}{\pi } \, k_{n} (x)}
  \\  \phantom{\rule{1pt}{15pt}}
   {\ds 4 x\left[k_{n} (x) h'_{ n\, -\, 2} (x)+k'_{ n\, -\, 2} (x) h_{n} (x)\right]=}
    \\  \phantom{\rule{1pt}{15pt}}
    {\ds (n-2)\, \left[k_{n} (x) h_{ n\, -\, 2} (x)+k_{ n\, -\, 2} (x) h_{n} (x)\right]+} 
    \\ \phantom{\rule{1pt}{15pt}}
    {(n+2)\, \left[k_{n} (x) h_{ n\, +\, 2} (x)+k_{ n\, +\, 2} (x) h_{n} (x)\right]}
     \\  \phantom{\rule{1pt}{15pt}}
     {\ds \left[k_{n} (x) h''_{ n} (x)-k''_{n} (x) h_{ n} (x)\right]=-\frac{2}{\pi  x} \, k_{n} (x) }, \end{array} 
\end{equation} 
where $n$  is an even integer.

 If we  consider the Havelock function in the special case
\begin{equation} \label{GrindEQ__59_} 
h_{n} (n x)=\frac{2}{\pi } \, \int _{0}^{\pi /2}\sin [n (x\tan \theta - \theta )]\, d\theta  =\frac{2}{\pi } \, \int _{0}^{\pi /2}\sin [n \alpha ]\, d\theta,   
\end{equation} 
then we recognize that the  sums of series of the Havelock  can be expressed by finite trigonometric integrals. 

\noindent For example from [42]
\begin{equation} \label{GrindEQ__60_} 
{\ds \frac{2}{\pi } \, \sum _{n\, =\, 1}^{\infty }t^{n} \sin (n \alpha )=\frac{2}{\pi } \, \left[\frac{t\sin \alpha }{1-2 t\cos \alpha +t^{2} } \right]} \, ;\quad t^{2} <1 
\end{equation} 
and integrating \eqref{GrindEQ__60_} with interchanging the order of summation and integration, we have 
\begin{equation} \label{GrindEQ__61_} 
{\ds \sum _{n\, =\, 1}^{\infty }t^{n} h_{n} (n x)=\frac{2}{\pi } \, \int _{0}^{\pi /2}\frac{t\sin (x\tan \theta -\theta )}{1-2 t\cos (x\tan \theta -\theta )+t^{2} }   \, d\theta}
\,  ;\quad t^{2} <1.
\end{equation} 
In  a
similar way it is possible to obtain for series of the Bateman functions
\begin{equation} \label{GrindEQ__62_} 
{\ds \sum _{n\, =\, 1}^{\infty }t^{n} k_{n} (n x)=\frac{2}{\pi } \, \int _{0}^{\pi /2}\frac{1-t\cos (x\tan \theta -\theta )}{1-2 t\cos (x\tan \theta -\theta )+t^{2} }   \, d\theta}
 \, ;\quad t^{2} <1. 
\end{equation} 
By this procedure, using various finite and infinite trigonometric series from [42], many sums of the Bateman $k_n(nx)$  and Havelock  $h_n(nx)$  series with different coefficients, can be expressed by corresponding integrals.

 \vs

\section{The Generalized Bateman and Havelock Functions with Integer Orders}

\noindent In order to solve dual, triple or multi series equations, a number of generalized Bateman and Havelock functions were introduced [25,26,29-35]. 
From the generalized functions only two considered in 1972 by Srivastava [31] are presented here. There is no agreed uniform notation of the generalized Bateman and Havelock functions. 
They are defined by using different letters, with upper and lower indexes. Here these functions are presented with an additional lower index 
with $k>-1$ as
\begin{equation} \label{GrindEQ__63_} 
\begin{array}{l} 
{\ds k_{n,k} (x)=\frac{2}{\pi } \, \int _{0}^{\pi /2}(\cos \theta )^{k} \cos (x\tan \theta -n \theta )\, d\theta, }
 \\  \phantom{\rule{1pt}{15pt}}
  {\ds h_{n,k} (x)=\frac{2}{\pi } \, \int _{0}^{\pi /2}(\cos \theta )^{k} \sin (x\tan \theta -n \theta )\, d\theta.  } 
  \end{array} 
\end{equation} 
It is suggested that if powers of cosine and sine functions appear also in \eqref{GrindEQ__63_}, then the third lower index ${m}$  is included
\begin{equation} \label{GrindEQ__64_} 
\begin{array}{l} 
{\ds k_{n,k,m} (x)=\frac{2}{\pi } \, \int _{0}^{\pi /2}(\cos \theta )^{k} (\sin \theta )^{m} \cos (x\tan \theta -n \theta )\, d\theta,  } 
\\  \phantom{\rule{1pt}{15pt}} 
{\ds h_{n,k,m} (x)=\frac{2}{\pi } \, \int _{0}^{\pi /2}(\cos \theta )^{k} (\sin \theta )^{m} \sin (x\tan \theta -n \theta )\, d\theta,  } 
\end{array} 
\end{equation} 
where this notation differs from that used in \eqref{GrindEQ__7_}.

 Values of three such integrals having $n = 0$ and $k = 0,1, 2$  are known
\begin{equation} \label{GrindEQ__65_} 
\begin{array}{l}
 {\ds \int _{0}^{\pi /2}(\cos \theta )^{2} \cos (x\tan \theta -n \theta )\, d\theta  =\frac{\pi (1+x)\, e^{-\, x} }{4} =\frac{\pi }{2} \, k_{0,2} (x)} 
 \\  \phantom{\rule{1pt}{15pt}}
  {\ds \int _{0}^{\pi /2}(\sin \theta )^{2} \cos (x\tan \theta -n \theta )\, d\theta  =\frac{\pi (1-x)\, e^{-\, x} }{4} =\frac{\pi }{2} \, k_{0,0,2} (x)}
   \\  \phantom{\rule{1pt}{15pt}}
    {\ds \int _{0}^{\pi /2}\cos \theta \sin \theta \sin (x\tan \theta -n \theta )\, d\theta  =\frac{\pi x\, e^{-\, x} }{4} =\frac{\pi }{2} \, h_{0,1,1} (x)}.
     \end{array} 
\end{equation} 

 The recurrence and differential expressions for the generalized Havelock  functions are [31]
\begin{equation} \label{GrindEQ__66_} 
\begin{array}{l}
 {\ds \left[(n-k-2)\, h_{ n\, -\, 2,k} (x)+(n+k+2)\, h_{ n\, +\, 2,k} (x)+(2 n-x)\, h_{n,k} (x)\right]=-\frac{8}{\pi } }
  \\ \phantom{\rule{1pt}{15pt}}
   {\ds 4 x h'_{n,k} (x)=\left[(n-k-2)\, h_{ n\, -\, 2,k} (x)-(n+k+2)\, h_{ n\, +\, 2,k} (x)+2 k\, h_{n,k} (x)\right]} 
   \\  \phantom{\rule{1pt}{15pt}}
   {\ds 2 x h'_{n,k} (x)-\frac{4}{\pi } =\left[(n-k-2)\, h_{ n\, -\, 2,k} (x)+(n+k-2 x)\, h_{ n\, +\, 2,k} (x)\right]}
    \\  \phantom{\rule{1pt}{15pt}}
    {\ds 2 h'_{0,2 k} (x)=\left[2 h_{0,2 k\, +\, 2} (x)-h_{0,2 k} (x)-h_{2,2 k\, +\, 2} (x)\right]}
     \\  \phantom{\rule{1pt}{15pt}}
      {\ds x h''_{n,k} (x)-k h'_{n,k} (x)+(n-x)h_{n,k} (x)=-\frac{2}{\pi } },
       \end{array} 
\end{equation} 
and for the generalized Bateman functions
\begin{equation} \label{GrindEQ__67_} 
2 k'_{0,2 k} (x)=\left[2 k_{0,2 k\, +\, 2} (x)-k_{0,2 k} (x)-k_{2,2 k\, +\, 2} (x)\right]\,. 
\end{equation} 
In 1972 Srivastava [31] was able to show that
\begin{equation} \label{GrindEQ__68_} 
\begin{array}{l} 
{\ds k_{0,2 k} (x)=\frac{2}{\pi } \, \int _{0}^{\pi /2}(\cos \theta )^{2 k} \cos (x\tan \theta )\, d\theta  =}
 \\  \phantom{\rule{1pt}{15pt}}
  {\ds \frac{2}{\sqrt{\pi } \, \Gamma (k+1)} \, \left(\frac{x}{2} \right)^{k\, +\, 1/2} K_{k\, +\, 1/2} (x)}
   \\  \phantom{\rule{1pt}{15pt}}
   {\ds h_{0,2 k} (x)=\frac{2}{\pi } \, \int _{0}^{\pi /2}(\cos \theta )^{2 k} \sin (x\tan \theta )\,
    d\theta  =}
     \\  \phantom{\rule{1pt}{15pt}}
      {\ds \frac{2 \Gamma (-k)}{\sqrt{\pi } } \, \left(\frac{x}{2} \right)^{k\, +\, 1/2} \left[I_{k\, +\, 1/2} (x)-L_{-\, k\, -\, 1/2} (x)\right]},
       \end{array} 
\end{equation} 
and in the explicit form for the generalized Havelock  function
\begin{equation} \label{GrindEQ__69_} 
h_{2 n,2 k} (x)=\frac{1}{\pi } \left[k_{2 n} (x)\, li(e^{x} )-2 S_{n\, -\, k\, -1,k} (x)\right]
\, ;\quad n\ge k+1, 
\end{equation} 
where he determined the following polynomials for the expression in \eqref{GrindEQ__69_}
\begin{equation} \label{GrindEQ__70_} 
\begin{array}{l} 
{\ds S_{2,1} (x)=\frac{1}{6} \left(2+x+x^{2} \right)} \\ {S_{3,1} (x)=\frac{1}{12} \left(2-x^{2} +x^{3} \right)}
 \\  \phantom{\rule{1pt}{15pt}}
 {\ds S_{4,1} (x)=\frac{1}{30} \left(4+x+2 x^{2} -4 x^{3} +x^{4} \right)}
  \\  \phantom{\rule{1pt}{15pt}}
   {\ds S_{5,1} (x)=\frac{1}{180} \left(18-9 x^{2} +31 x^{3} -16 x^{4} +2 x^{5} \right)}
    \\  \phantom{\rule{1pt}{15pt}} 
    {\ds S_{5,1} (x)=\frac{1}{180} \left(18-9 x^{2} +31 x^{3} -16 x^{4} +2 x^{5} \right)}, \end{array} 
\end{equation} 
and
\begin{equation} \label{GrindEQ__71_} 
\begin{array}{l} 
{\ds S_{3,2} (x)=\frac{1}{48} \left(16+7 x+3 x^{2} +x^{3} \right)}
 \\  \phantom{\rule{1pt}{15pt}}
  {\ds S_{4,2} (x)=\frac{1}{120} \left(24+6 x+2 x^{2} +x^{3} +x^{4} \right)} 
  \\  \phantom{\rule{1pt}{15pt}} 
  {\ds S_{5,2} (x)=\frac{1}{360} \left(48+6 x- x^{3} -2 x^{4} +x^{5} \right)} \\ {S_{6,2} (x)=\frac{1}{2520} \left(268+30x+6 x^{2} +5 x^{3} +11 x^{4} -44 x^{5} +2 x^{6} \right)}. \end{array} 
\end{equation} 
Besides, in 1972 H.M. Srivastava [31] evaluated four Laplace transforms of the generalized Bateman and Havelock functions. Two are presented here, long but  complex expressions for the functions $k_{2,2k}(x)$ and $h_{2,2k}(x)$  are omitted here: 
\begin{equation} \label{GrindEQ__72_} 
\begin{array}{l} 
{\ds L\left\{k_{0,2 k} (x)\right\}=\left[\frac{(1-s)}{(1-s^{2} )^{k\, +\, 1} } -\frac{s}{\sqrt{\pi } } \, \sum _{m\, =\, 1}^{k}\frac{\Gamma (k-m+3/2)}{\Gamma (k-m+2)\, (1-s^{2} )^{m} }  \right]}
 \\ \phantom{\rule{1pt}{15pt}}
  {\ds L\left\{h_{0,2 k} (x)\right\}=\frac{1}{\pi } \left[\frac{2 \ln (s)}{(1-s^{2} )^{k\, +\, 1} } +\sum _{m\, =\, 1}^{k}\frac{1}{(k-m+1)\, (1-s^{2} )^{m} }  \right]}.
   \end{array} 
\end{equation} 
 For the solution of pairs of dual equations, other researchers called Srivastava [28,29] reported 
 a few more properties of the generalized Bateman functions, but these functions are slightly modified in their definitions.  

\vs  

\section {The Bateman and Havelock Functions with Unrestricted Orders}

\noindent General case of the Bateman and Havelock with any order 
\begin{equation} \label{GrindEQ__73_} 
\begin{array}{l}
 {\ds k_{\nu } (x)=\frac{2}{\pi } \, \int _{0}^{\pi /2}\cos (x\tan \theta -\nu  \theta )\, d\theta  }
  \\ \phantom{\rule{1pt}{15pt}}
   {\ds h_{\nu } (x)=\frac{2}{\pi } \, \int _{0}^{\pi /2}\sin (x\tan \theta -\nu  \theta )\, d\theta  } \end{array} 
\end{equation} 
is practically unknown in the literature, with only one exception, the definition of the Bateman function in 
 terms of the Whittaker function $W_{k,\mu}(z)$  or Tricomi function
 $U(a,b,z)$  (particular cases of the confluent hypergeometric function ) [7]
\begin{equation} \label{GrindEQ__74_} 
\begin{array}{l}
 {\ds k_{2 \nu } (x)=\frac{1}{\Gamma (\nu +1)} \, W_{\nu ,1/2} (2 x)=\frac{e^{-\, x} }{\Gamma (\nu +1)} \, U(-\, \nu ,0;2 x)}
  \\ \phantom{\rule{1pt}{15pt}}
   {\ds U(-\, \nu ,0;2 x)=2 x\, U(1-\nu ,2;2 x)} 
   \\  \phantom{\rule{1pt}{15pt}}
   {\ds k_{2 n\, +\, 2} (x)=2 x e^{-\, x} {}_{1} F_{1} (-2n;2;2x)\, ;\quad n=0,1,2,3,...} \end{array} 
\end{equation} 
Evidently, the corresponding generalized functions are
\begin{equation} \label{GrindEQ__75_} 
\begin{array}{l}
 {\ds k_{\nu ,\alpha ,\beta } (x)=\frac{2}{\pi } \, \int _{0}^{\pi /2}(\cos \theta )^{\alpha } (\sin \theta )^{\beta } \cos (x\tan \theta -\nu  \theta )\, d\theta  }
  \\  \phantom{\rule{1pt}{15pt}}
  {\ds h_{\nu ,\alpha ,\beta } (x)=\frac{2}{\pi } \, \int _{0}^{\pi /2}(\cos \theta )^{\alpha } (\sin \theta )^{\beta } \sin (x\tan \theta -\nu  \theta )\, d\theta,  } 
  \end{array} 
\end{equation} 
where $\alpha$, $\beta$ and $\nu$   have any real value. 
By changing the integration variable in \eqref{GrindEQ__73_} and \eqref{GrindEQ__75_}, 
$t = tan(\theta)$  , these functions can be expressed by infinite integrals
\begin{equation} \label{GrindEQ__76_} 
\begin{array}{l}
 {\ds k_{\nu } (x)=\frac{2}{\pi } \, \int _{0}^{\infty }\frac{\left[\cos (xt)\, \cos [\nu \tan ^{-\, 1} (t)]+\sin (xt)\, \sin [\nu \tan ^{-\, 1} (t)]\right]}{1+t^{2} }  \, dt}
  \\ \phantom{\rule{1pt}{15pt}}
   {\ds k_{\nu ,\alpha ,\beta } (x)=\frac{2}{\pi } \, \int _{0}^{\infty }\frac{t^{\beta } \left[\cos (xt)\, \cos [\nu \tan ^{-\, 1} (t)]+\sin (xt)\, \sin [\nu \tan ^{-\, 1} (t)]\right]}{(1+t^{2} )^{\alpha /2\, +\, \beta /2\, +1} }  \, dt} 
   \\ \phantom{\rule{1pt}{15pt}}
   {\ds h_{\nu } (x)=\frac{2}{\pi } \, \int _{0}^{\infty }\frac{\left[\sin (xt)\, \cos [\nu \tan ^{-\, 1} (t)]-\cos (xt)\, \sin [\nu \tan ^{-\, 1} (t)]\right]}{1+t^{2} }  \, dt}
    \\  \phantom{\rule{1pt}{15pt}}
     {\ds h_{\nu ,\alpha ,\beta } (x)=\frac{2}{\pi } \, \int _{0}^{\infty }\frac{t^{\beta } \left[\sin (xt)\, \cos [\nu \tan ^{-\, 1} (t)]-\cos (xt)\, \sin [\nu \tan ^{-\, 1} (t)]\right]}{(1+t^{2} )^{\alpha /2\, +\, \beta /2\, +1} }  \, dt.}
      \end{array} 
\end{equation} 

 In Figures 5 and 6  
 we  illustrate the  behavior and the symmetries with respect to the order of the Bateman functions with fractional positive and negative values: 
 $\nu = n+1/2$  and $\nu = - ({n} + 1/2) $,  with ${n} = 0, 1, 2, 3, 4, 5$. 
 The same is demonstrated in Figures 7 and 8 for the Havelock-functions. 
 Similarly as in \eqref{GrindEQ__28_}, differentiation of the Bateman functions 
 with respect to the argument $x$  is for $k=0,1,2,3,...$

\begin{equation} \label{GrindEQ__77_} 
\begin{array}{l} 
{\ds \frac{\partial ^{2 k} \, k_{\nu } (x)}{\partial  x^{2 k} } 
=
(-1)^{k} \frac{2}{\pi } \, \int _{0}^{\pi /2}(\tan \theta )^{2 k} 
\cos (x\tan \theta -\nu \theta )\, d\theta  }
 \\  \phantom{\rule{1pt}{15pt}}
  {\ds \frac{\partial ^{2 k\, +\, 1} \, k_{\nu } (x)}{\partial  x^{2 k\, +\, 1} } =(-1)^{k} \frac{2}{\pi } \, \int _{0}^{\pi /2}(\tan \theta )^{2 k\, +\, 1} \sin (x\tan \theta -\nu \theta )\, d\theta.  }
  \end{array} 
\end{equation}

\begin{figure}[h!]
	\centering
	\includegraphics[width=12cm]{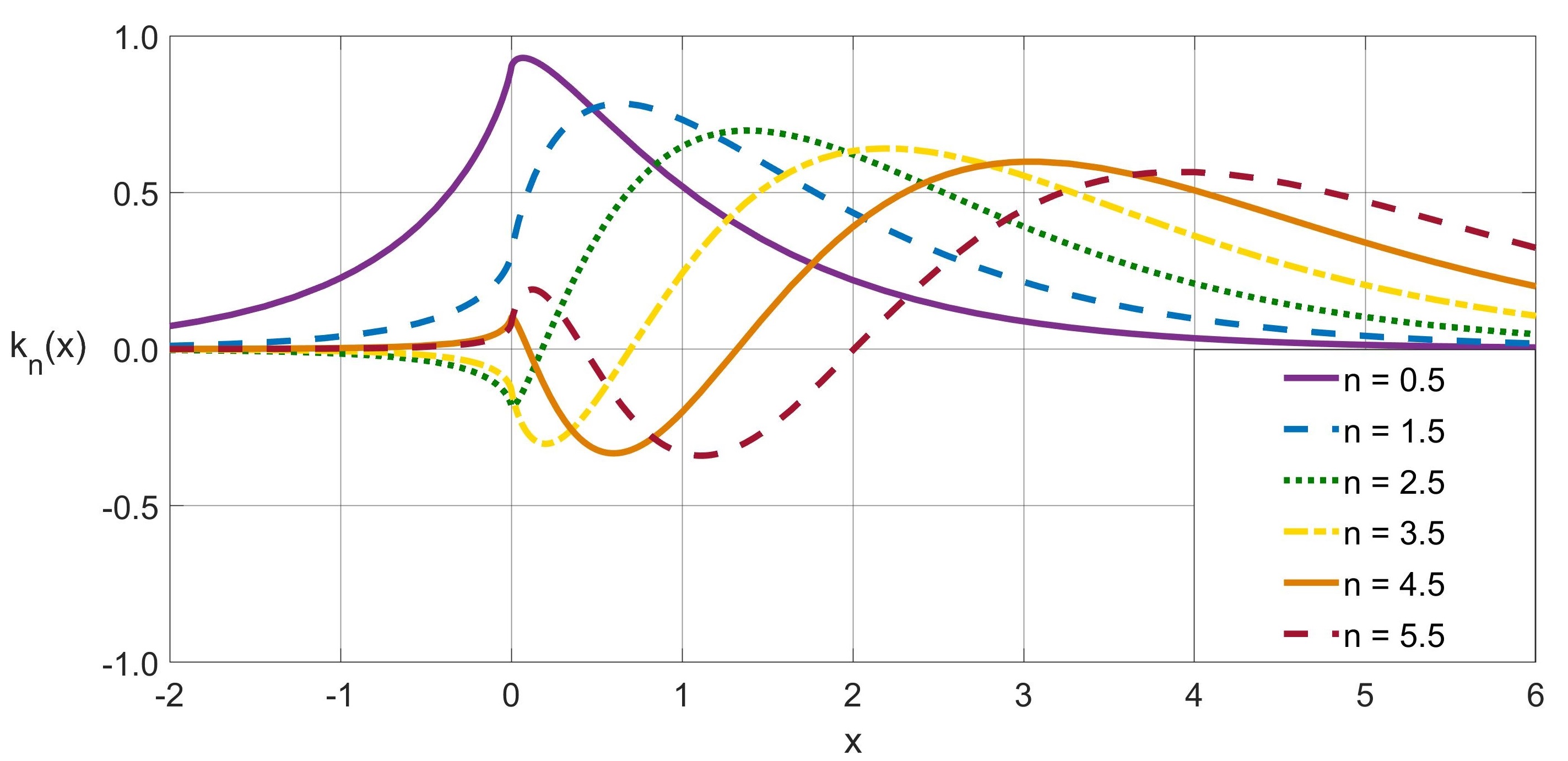}
	\caption{Bateman functions with  positive \textit{n} + 1/2 orders as a function of argument \textit{x}.} 
\end{figure}

\begin{figure}[h!]
	\centering
	\includegraphics[width=12cm]{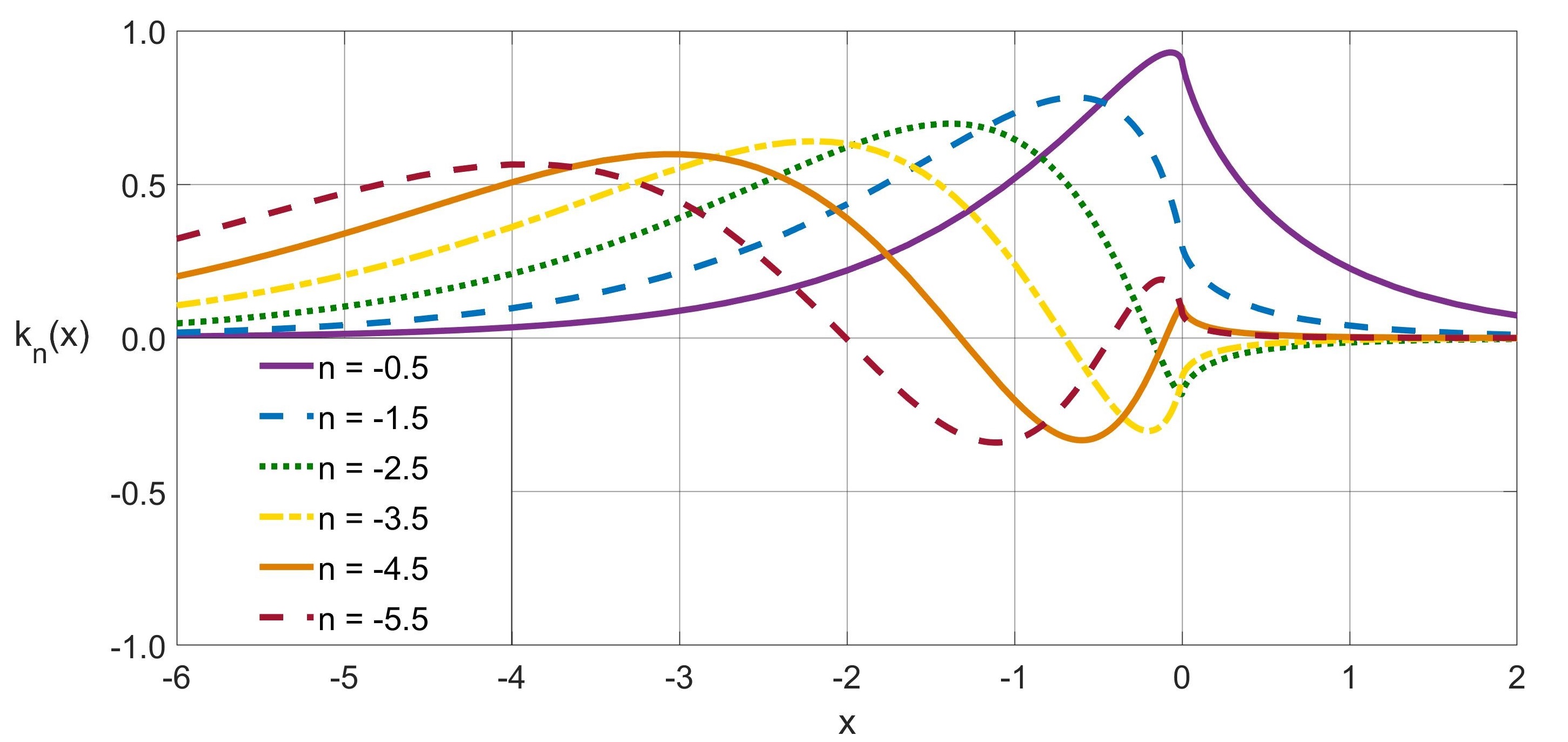}
	\caption{Bateman functions with  negative \textit{n} + 1/2 orders as a function of argument \textit{x}.} 
\end{figure}

\begin{figure}[h!]
	\centering
	\includegraphics[width=12cm]{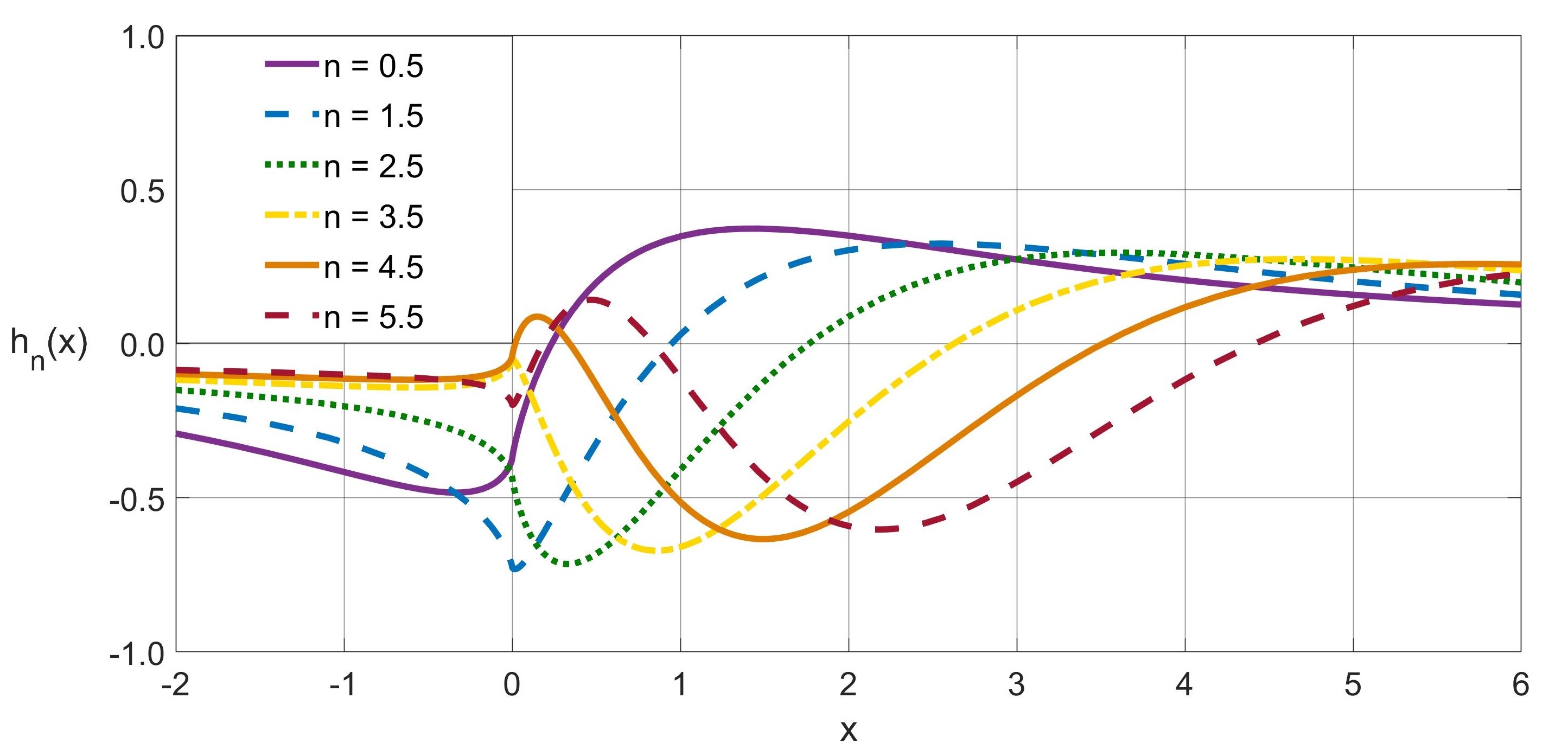}
	\caption{ Havelock  functions with  positive \textit{n} + 1/2 orders as a function of argument \textit{x}.} 
\end{figure}

\begin{figure}[h!]
	\centering
	\includegraphics[width=12cm]{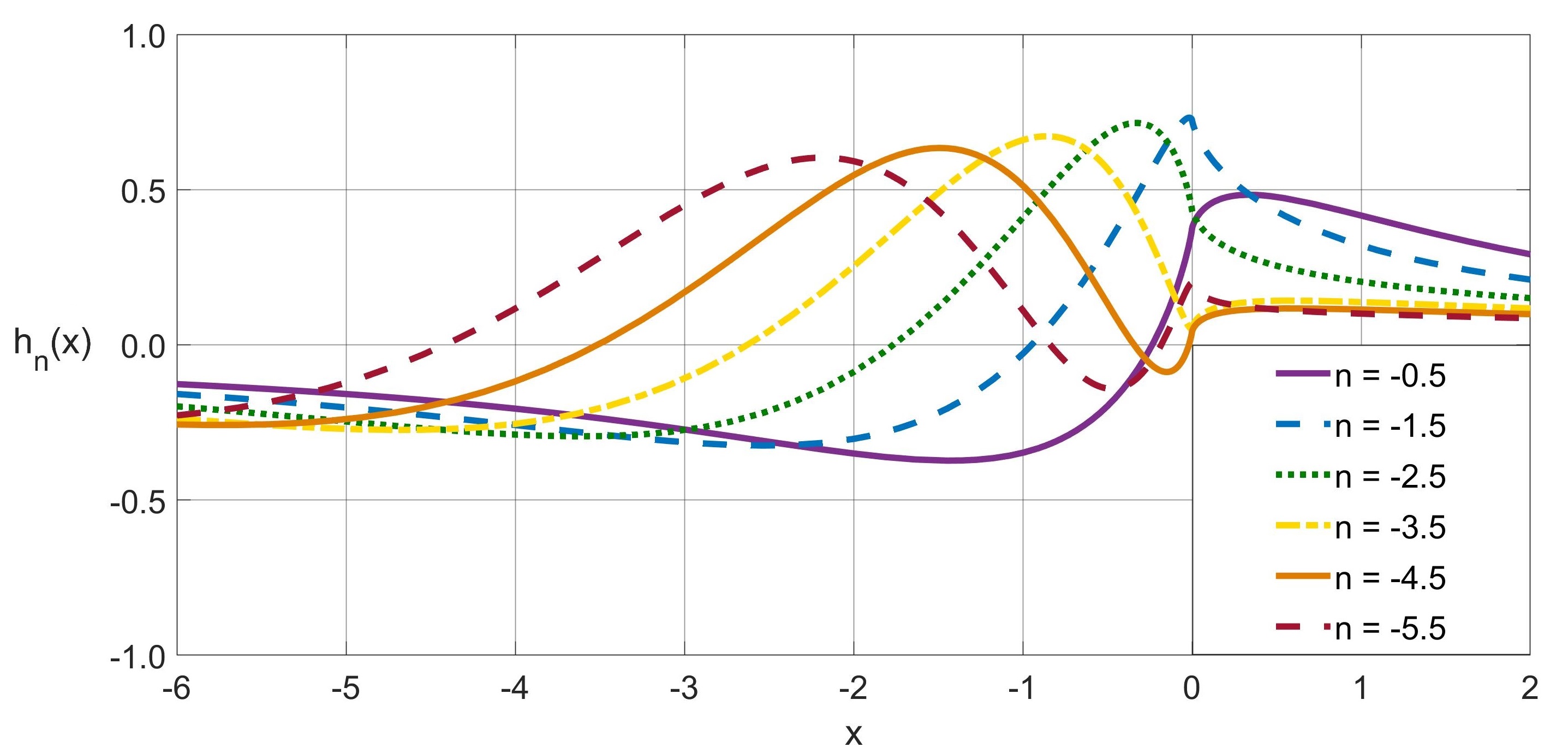}
	\caption{Havelock  functions with  negative \textit{n} + 1/2 orders as a function of argument \textit{x}.} 
\end{figure}
\vskip 0.5truecm
\noindent and in the case of the Havelock  functions 

\begin{equation} \label{GrindEQ__78_} 
\begin{array}{l} 
{\ds \frac{\partial ^{2 k} \, h_{\nu } (x)}{\partial  x^{2 k} } =(-1)^{k} \frac{2}{\pi } \,
 \int _{0}^{\pi /2}(\tan \theta )^{2 k} \sin (x\tan \theta -\nu  \theta )\, d\theta  }
  \\  \phantom{\rule{1pt}{15pt}}
   {\ds \frac{\partial ^{2 k\, +\, 1} \, h_{\nu } (x)}{\partial  x^{2 k\, +\, 1} } =(-1)^{k} \frac{2}{\pi } \, \int _{0}^{\pi /2}(\tan \theta )^{2 k\, +\, 1} \cos (x\tan \theta -\nu  \theta )\, d\theta  }
    \\  \phantom{\rule{1pt}{15pt}}
    {\ds k=0,1,2,3,...}
     \end{array} 
\end{equation} 
 Using the 
  definition of these function from \eqref{GrindEQ__73_}, it is possible to consider the Bateman and Havelock functions as functions of two variables 
 ${x}$ and ${\nu}$. 
 Thus, it is possible also to perform differentiation with respect to ${\nu}$
\begin{equation} \label{GrindEQ__79_} 
 \begin{array}{l} 
 {\ds \frac{\partial ^{2 k} \, k_{\nu } (x)}{\partial  \nu ^{2 k} } =(-1)^{k} \frac{2}{\pi } \,
  \int _{0}^{\pi /2}\theta ^{2 k} \cos (x\tan \theta -\nu  \theta )\, d\theta  }
  \\  \phantom{\rule{1pt}{15pt}} 
  {\ds \frac{\partial ^{2 k\, +\, 1} \, k_{\nu } (x)}{\partial  \nu ^{2 k\, +\, 1} } =(-1)^{k} \frac{2}{\pi } \, \int _{0}^{\pi /2}\theta ^{2 k\, +\, 1} \sin (x\tan \theta -\nu  \theta )\, d\theta  }
   \\  \phantom{\rule{1pt}{15pt}}
   {\ds k=0,1,2,3,...}
    \end{array}
   \end{equation}

\noindent and
\begin{equation} \label{GrindEQ__80_} 
 \begin{array}{l} 
 {\ds \frac{\partial ^{2 k} \, h_{\nu } (x)}{\partial  \nu ^{2 k} } =(-1)^{k} \frac{2}{\pi } \,
  \int _{0}^{\pi /2}\theta ^{2 k} \sin (x\tan \theta -\nu  \theta )\, d\theta  }
  \\  \phantom{\rule{1pt}{15pt}}
   {\ds \frac{\partial ^{2 k\, +\, 1} \, h_{\nu } (x)}{\partial  \nu ^{2 k\, +\, 1} } =(-1)^{k} \frac{2}{\pi } \, \int _{0}^{\pi /2}\theta ^{2 k\, +\, 1}
    \cos (x\tan \theta -\nu  \theta )\, d\theta  }
     \\  \phantom{\rule{1pt}{15pt}}
     {\ds k=0,1,2,3,...}
      \end{array}
     \end{equation}

\noindent 
the
first derivatives with respect to the order at fixed positive and negative values of argument ${x}$ of the Bateman functions are plotted in Figure 9 and 10, and the same for the Havelock functions in Figure 11 and 12. 
As can be observed, these functions are symmetrical in both cases.

 If orders are pure imaginary numbers $\nu = i\alpha$  then the Bateman and Havelock  functions become complex functions which are expressed by integrals with integrands having products of trigonometric and hyperbolic functions.
 
 As pointed out above, the Bateman and Havelock functions were introduced to the mathematical literature as solutions of particular problems in fluid mechanics [20,36]. 

\begin{figure}[h!]
	\centering
	\includegraphics[width=12cm]{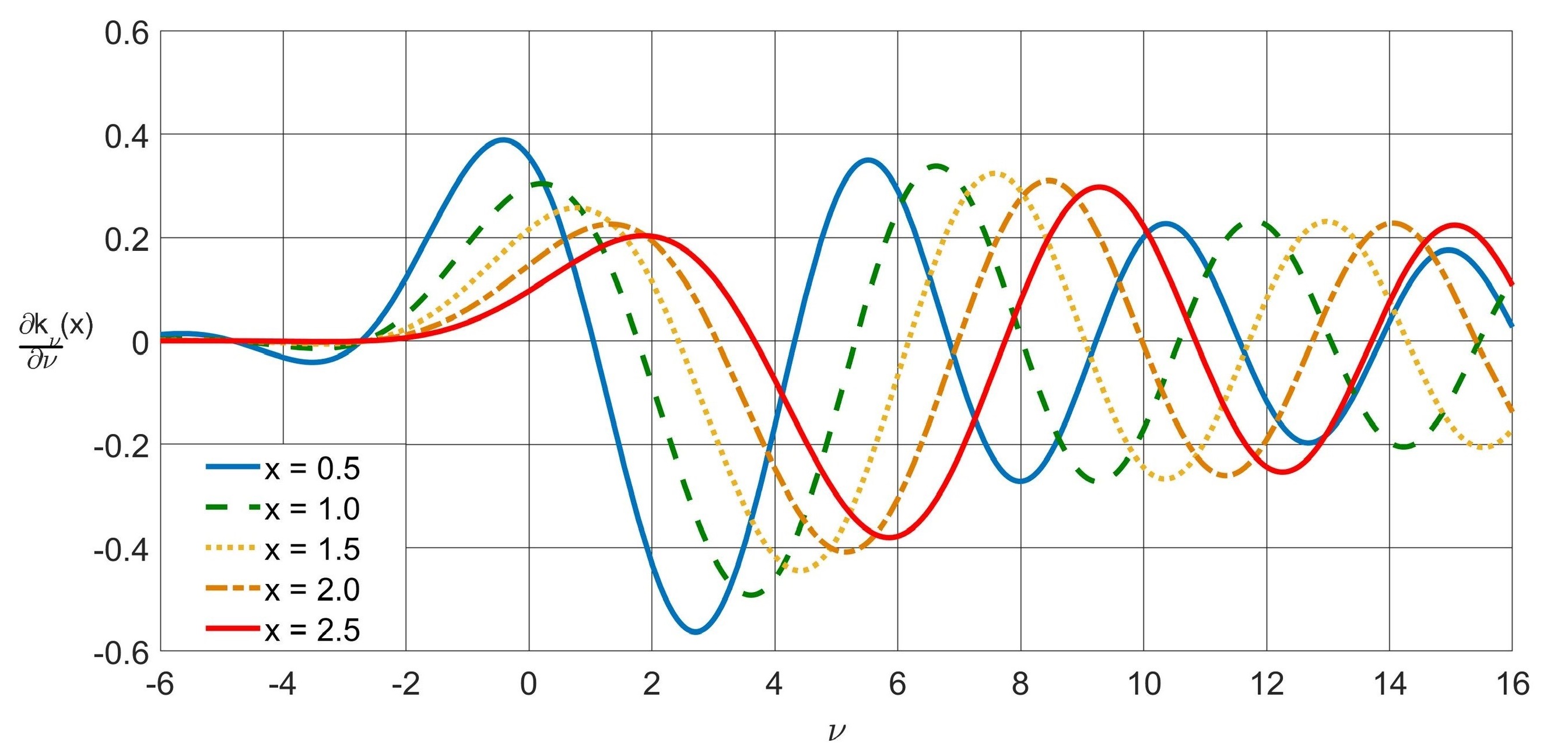}
	\caption{First derivatives of the Bateman functions with respect to the order at fixed positive values of argument \textit{x}.} 
\end{figure}

\begin{figure}[h!]
	\centering
	\includegraphics[width=12cm]{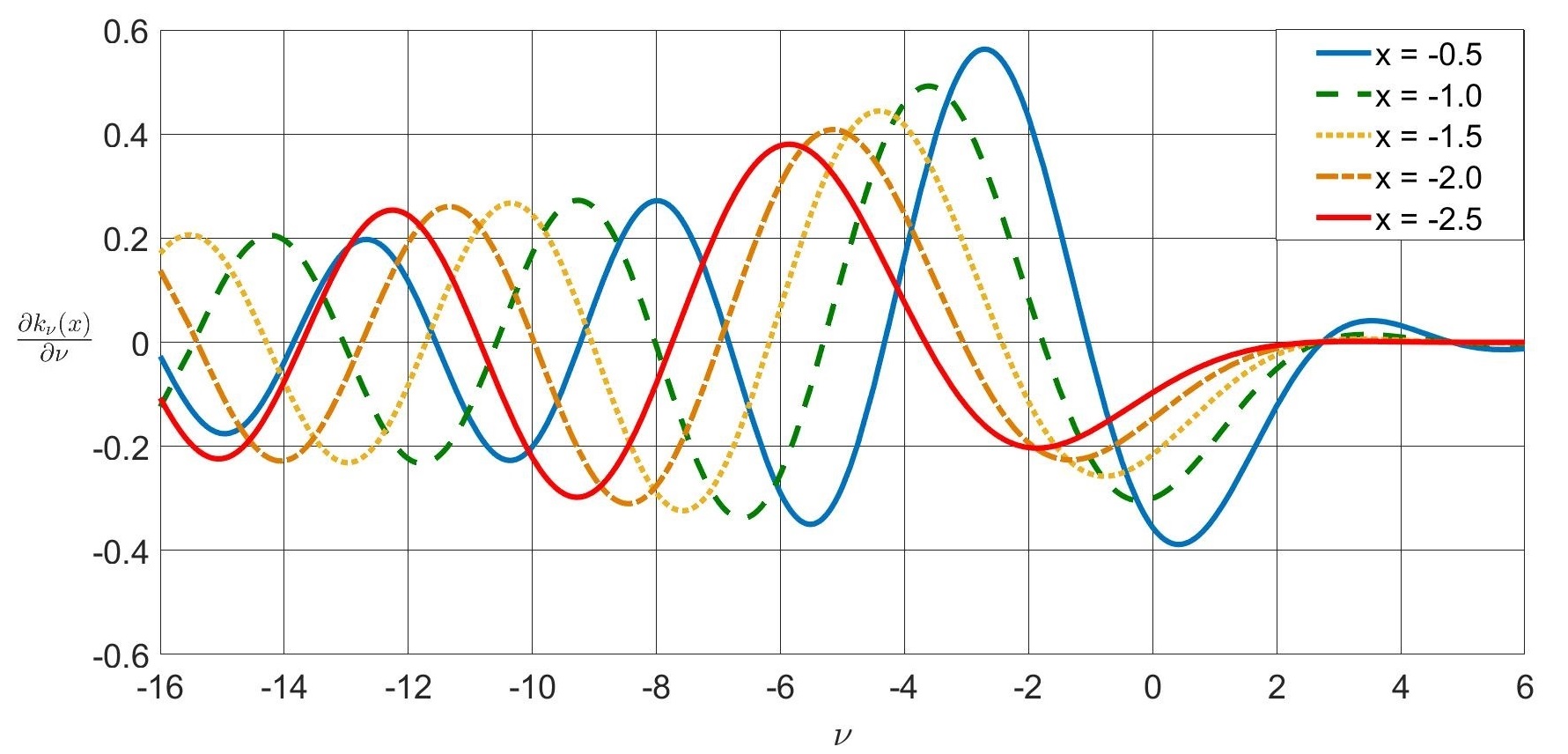}
	\caption{First derivatives of the Bateman functions with respect to the order at fixed negative values of argument \textit{x}.} 
\end{figure}

\begin{figure}[h!]
	\centering
	\includegraphics[width=12cm]{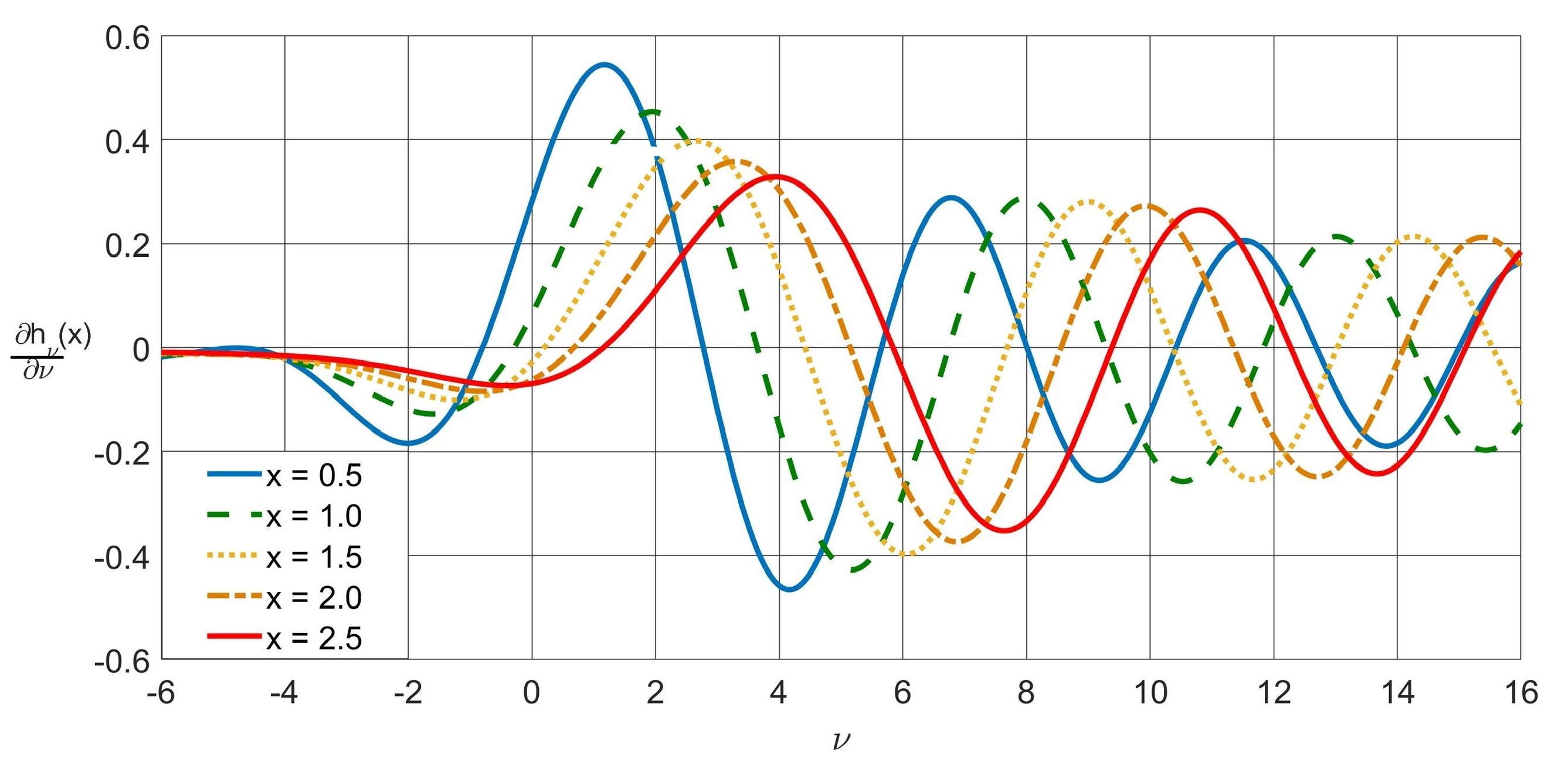}
	\caption{First derivatives of the Havelock  functions with respect to the order at fixed positive values of argument \textit{x}.} 
\end{figure}

\begin{figure}[h!]
	\centering
	\includegraphics[width=12cm]{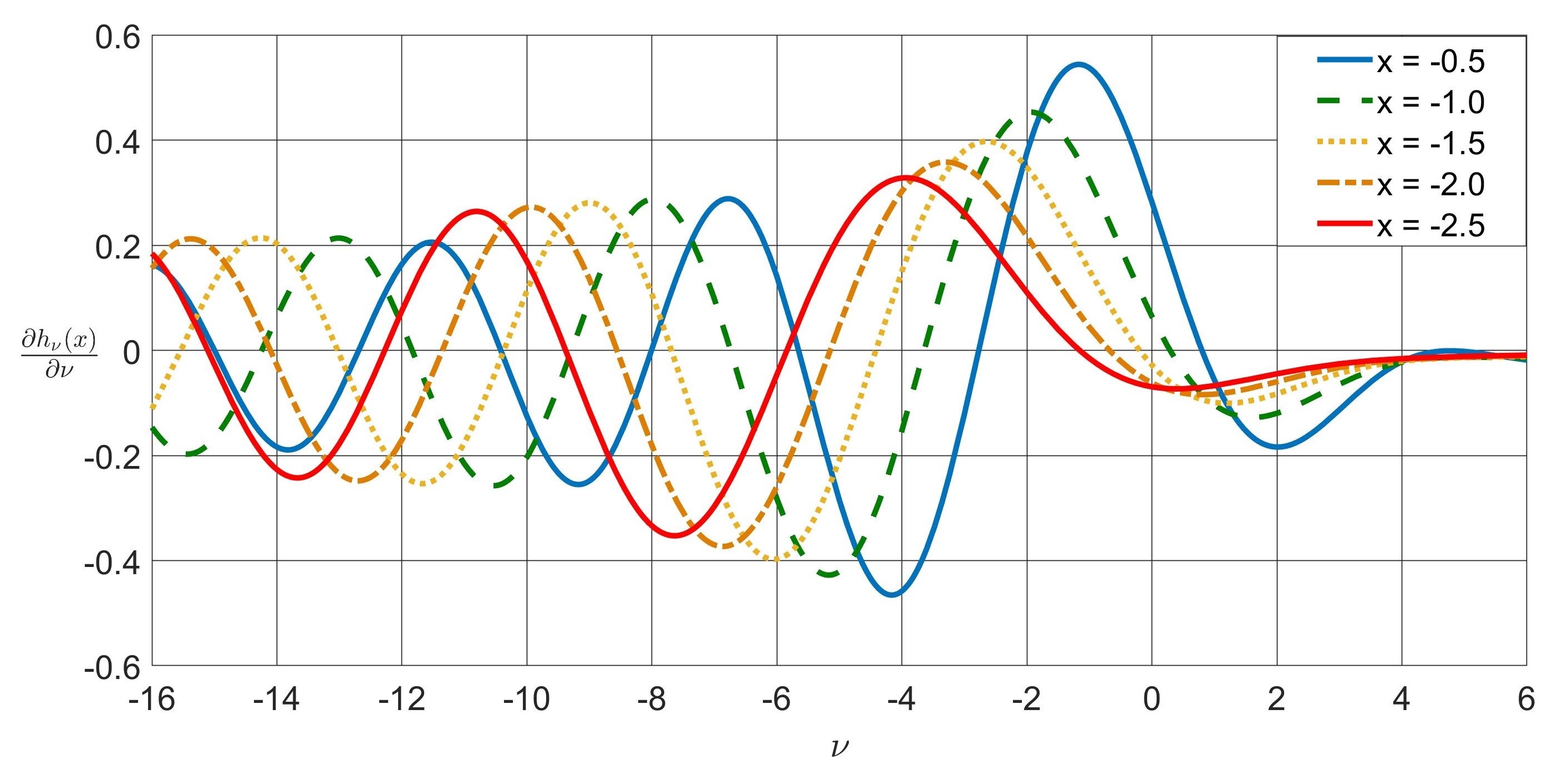}
	\caption{First derivatives of the Havelock functions with respect to the order at fixed negative values of argument \textit{x}.} 
\end{figure}

Years later, these functions were generalized to the form given in 
  \eqref{GrindEQ__64_}
  and \eqref {GrindEQ__75_}
 [25,26,29-35]. 
 It should be mentioned however, that historically, these proposed generalizations are not new, and they were already discussed much earlier by Giuliani in1888 [43] and by Bateman in 1931 [44]. 
 They also introduced similar trigonometric integrals, but in the context of particular cases of the Kummer confluent hypergeometric functions. 
 It is rather strange, that in the later investigations [25,26,29-35], when the generalized Bateman and Havelock functions were proposed, previous studies on this subject were completely ignored. 
 Considering that the trigonometric integrals and associated with them differential equations 
presented in the Giuliani and Bateman papers  
 are of particular importance and interest, it was decided to summarize their results separately, in Appendix B.

\vs

\section{The Bateman-Integral Functions}

\noindent Analogous to the sine-integral, cosine integral and the Bessel-integral functions
\begin{equation} \label{GrindEQ__81_} 
\begin{array}{l}
 {\ds si(x)=-\, \int _{x}^{\infty }\frac{\sin t}{t}  \, dt} \\ {Ci(x)=-\, \int _{x}^{\infty }\frac{\cos t}{t}  \, dt}
  \\ \phantom{\rule{1pt}{15pt}}
  {\ds Ji_{\nu } (x)=-\, \int _{x}^{\infty }\frac{J_{\nu } (t)}{t}  \, dt.}
   \end{array} 
\end{equation} 
Chaudhuri [26] has  introduced the Bateman-integral function
\begin{equation} \label{GrindEQ__82_} 
 {\ds ki_{2 n} (x)=-\, \int _{x}^{\infty }\frac{k_{2 n} (t)}{t}  \, dt\, ;\quad t>0,}   
\end{equation}
 and mainly using operational calculus 
 he has 
 discussed its properties. 
\begin{equation} \label{GrindEQ__83_} 
 \begin{array}{l} 
 {\ds ki_{2 n} (x)=\, \int _{0}^{x}\frac{k_{2 n} (t)}{t}  \, dt+ki_{2 n} (0)} 
 \\ \phantom{\rule{1pt}{15pt}}
  {\ds ki_{2 n} (0)=0\quad ;\quad n=2k\,;\quad k=0,1,2,3,...}
   \\  \phantom{\rule{1pt}{15pt}}
   {\ds ki_{2 n} (0)=-\frac{2}{n} \, ;\quad n=2k+1}
    \end{array}
    \end{equation}
Using similarity  with  the Laguerre polynomials, Chaudhuri [26] derived the following series expressions for the Bateman-integral functions 
\begin{equation} \label{GrindEQ__84_} 
\begin{array}{l} 
{\ds ki_{2 n} (x)=\frac{e^{-\, x} }{n} \, \sum _{k\, =\, 1}^{n}(-2)^{k}  
\left(\begin{array}{c} {n} \\ {k} \end{array}\right)\, L_{k\, -\, 1} (x)}
 \\  \phantom{\rule{1pt}{15pt}}
  {\ds ki_{2 n} (x)=} 
  \\  \phantom{\rule{1pt}{15pt}}
   {\ds \frac{1}{n x} \, \left[n\, k_{2 n} (x)-2 \sum _{m\, =\, 1}^{n}(-1)^{k} \left(\begin{array}{c} {n}
    \\
     {\ds m} \end{array}\right)[m\, k_{2 m} (2 x)+(m+1) k_{2 m\, +\, 2} (2 x)-2k_{0} (2 x) \right]} 
     \\ \phantom{\rule{1pt}{15pt}}
      {\ds ki_{2 n} (x)=\frac{(-1)^{n\, -\, 1} e^{x} }{2^{n\, +\, 1} }
       \left[\sum _{m\, =\, 1}^{n}m\, k_{2 k} (x) \right]}
        \\  \phantom{\rule{1pt}{15pt}}
         {\ds L_{n\, -\, 1} (x)=\frac{e^{x} }{2^{n} } \, \sum _{m\, =\, 1}^{n}(-1)^{m}  \left(\begin{array}{c} {n} \\ {m} \end{array}\right)\, m\, ki_{2 m} (x)} \\ {ki_{2} (x)=-2\, k_{0} (x),}
          \end{array} 
\end{equation} 
and the recurrence and differential expressions
\begin{equation} \label{GrindEQ__85_} 
\begin{array}{l}
 {\ds k_{2 n} (x)=\frac{(n-1) ki_{2 n\, -\, 2} (x)-(n+1) ki_{2 n\, +\, 2} (x)}{2} }
  \\  \phantom{\rule{1pt}{15pt}}
   {\ds n\, ki_{2 n} (x)\, +\, (n+1)\, ki_{2 n\, +\, 2} (x)=-2\, \sum _{k\, =\, 0}^{n}ki_{2 k} (x) }
    \\  \phantom{\rule{1pt}{15pt}}
     {\ds x\, ki'_{2 n} (x)=\frac{(n-1) ki_{2 n\, -\, 2} (x)-(n+1) ki_{2 n\, +\, 2} (x)}{2} }
      \\  \phantom{\rule{1pt}{15pt}}
       {\ds x\, ki'_{2 n} (x)=k_{2 n} (x).}
        \end{array} 
\end{equation} 
He was also able to relate the Bateman-integral functions with the Bessel and
 the
 Bessel integral functions
\begin{equation} \label{GrindEQ__86_} 
\begin{array}{l}
 {\ds (n+1)\, \left[Ji_{n\, +\, 1} (x) ki_{2 n\, -\, 2} (x)-Ji_{n\, -\, 1} (x) ki_{2 n\, +\, 2} (x)\right]=} 
 \\  \phantom{\rule{1pt}{15pt}}
 {\ds 2 x Ji_{n\, -\, 1} (x)  ki'_{2 n} (x)-2 n Ji'_{n} (x) ki_{2 n\, -\, 2} (x)}
  \\ \phantom{\rule{1pt}{15pt}}
   {\ds \sum _{m\, =\, 1}^{\infty }(-1)^{m}  m\, ki_{2 m} (x)\, ki_{2 m} (y)
   =J_{0} (2 \sqrt{x y} )}. 
   \end{array} 
\end{equation} 
From integral expressions, the Laplace transform are presented here, when indefinite, definite and infinite integrals related to of the Bateman-integral functions 
are given in Appendix A:
\begin{equation} \label{GrindEQ__87_} 
\begin{array}{l} 
{\ds L\left\{ki_{2 n} (x)\right\}=\frac{1}{n s} \, \left[\left(\frac{1-s}{s+1} \right)^{n} -1\right]=\frac{1}{n s} \, \sum _{k\, =\, 1}^{n}(-1)^{k}  
\left(\begin{array}{c} {n}  \\ {k} \end{array}\right)\, 
\left(\frac{2 s}{s+1} \right)^{k} }
 \\  \phantom{\rule{1pt}{15pt}}
  {\ds L\left\{ki_{2 n} (2 x)\right\}=
  \frac{1}{n s} \, \left[\left(\frac{2-s}{s+2} \right)^{n} -1\right]} 
  \\ \phantom{\rule{1pt}{15pt}}
   {\ds L\left\{ki_{0} (x)\right\}=-\frac{\ln (s)}{s} } 
   \\  \phantom{\rule{1pt}{15pt}}
   {\ds L\left\{ki_{2} (x)\right\}=-\, \frac{2}{s+1}. }
    \end{array} 
\end{equation} 
It is also worthwhile to mention that Srivastava [25] expressed the Bateman-integral function in 
the
following way
\begin{equation} \label{GrindEQ__88_} 
\begin{array}{l} 
{\ds ki_{2 n} (x)=\frac{\pi }{2} \, \left[k'_{2 n} (x) h_{2 n} (x)-h'_{2 n} (x) k_{2 n} (x)\right]=}
 \\  \phantom{\rule{1pt}{15pt}}
  {\ds \frac{\pi }{8 x} \left[(2n+2) [k_{2 n} (x) h_{2 n\, +\, 2} (x)-k_{2 n\, +\, 2} (x)h_{2 n} (x)]\right. -}
   \\  \phantom{\rule{1pt}{15pt}}
    {\ds \left. (2n-2) [k_{2 n} (x) h_{2 n\, -\, 2} (x)-k_{2 n\, -\, 2} (x)h_{2 n} (x)]\right]}, \end{array} 
\end{equation} 
by including products of the Bateman and Havelock  functions.

\vs

\section {Conclusions}

\noindent As solutions of fluid mechanics problems, 
more than ninety years ago, Havelock in 1925 and Bateman in 1931 introduced  new functions which are expressed in terms of finite trigonometric integrals and discussed their properties. Initially, these functions found attention of a number of mathematicians who further developed this subject and proposed some generalizations.
   However,  unfortunately, after a rather short period, the Havelock  and Bateman functions were practically abandoned. 
Today, only the Bateman function is listed in mathematical handbooks as a particular case of the confluent hypergeometric function, thus as a minor special function. However, as is clearly showed in this survey, these functions have interesting properties and a rather large mathematical material was devoted and associated with them. 
This leads to conclusion that they should be treated as independent special functions. 
Since at present, in reference books, our knowledge about these functions is very limited,
 we decided to prepare this survey where basic properties of the Havelock and Bateman functions are presented.
  We have found useful for  the reader's convenience to add two Appendixes:
 Appendix A is devoted to integrals 
 associated with the Bateman and Bateman-integral functions
 whereas  Appendix B  is devoted to
 trigonometric integrals and differential equations associated with the Kummer Confluent Hypergeometric Functions according to the almost unknown papers by Giuliani [43[ and by Bateman himself [44]. 
\\
In Appendix C   we have added the integral representations of the special functions used in this survey.
\\
 It is worth to note that the {Bateman Manuscript} is currently under revision with 
  the name {\it Encyclopedia of Special Functions: the Askey-Bateman Project}, see [45].
  However the volume dealing with  the confluent hypergeometric functions is not yet available.
  \\ \\
 \vvs
\noindent
{\bf Funding}:
{This research received no external funding.} 
\\
{\bf Conflicts of interest}: {The authors declare no conflict of interest.}

\section*{Acknowledgments}The research  of FM and AC has been carried out in the framework of the activities of the National Group of Mathematical Physics (GNFM, INdAM).\\
All the the authors like to acknowledge the librarians of the Department of Physics and Astronomy of the University of Bologna to have found the pdf of several articles cited in the bibliography.
The reader is kindly requested to accept the authors'  somewhat informal style and 
to contact the corresponding author for pointing out possible misprints and mathematical errors.

\newpage

\section*{Appendix A: Integrals Associated with the Bateman and Bateman-Integral Functions}

\noindent 
The integrals presented here are compiled from the literature and they have a definite form. Their number can be enlarged by applying interconnections between the Bateman, Bateman-integral and other special functions and using operational calculus. Besides, there are many integrals which are expressed in term of infinite series, but they are omitted from this tabulation. 

$$
\int _{0}^{1}(1-t)^{\beta \, -\, 1}  \, e^{\alpha  t} \, k_{2 n} (\alpha  t)\, dt=\frac{(-1)^{n\, -\, 1} (n-1)!\, \Gamma (\beta )}{\Gamma (\beta +n+1)} \, L_{n\, -\, 1}^{(\beta \, +\, 1)} (2 \alpha )\, ;\quad \beta >0          
  \eqno (A.1)
 $$

 $$ \begin{array}{l} 
 {\ds \int _{0}^{x}k_{2 m}  (t)\, k_{2 n} (x-t)\, dt}=
 {\ds \int _{0}^{x}k_{2 n}  (t)\, k_{2 m} (x-t)\, dt}=
 \\  \phantom{\rule{1pt}{15pt}}
 {\ds \frac{1}{2} \, \left[k_{2 m\, +\, 2 n\, -2} (x)+2 k_{2 m\, +\, 2 n} (x)+k_{2 m\, +\, 2 n\, +2} (x)\right]} 
 \end{array}  
     \eqno (A.2)
 $$

 $$
 \begin{array}{l}
 {\ds \int _{0}^{x}\frac{J_{0} (t)-k_{0} (t)}{t}  \, dt}
 ={\ds Ji_{0} (x)-ki_{0} (x)+\ln 2}
 \\       \phantom{\rule{1pt}{15pt}}            
{\ds \int _{0}^{x}\frac{J_{n} (t)-k_{2 n} (t)}{t}  \, dt}
= {\ds Ji_{n} (x)-ki_{2 n} (x)+\frac{(-1)^{n} }{n}} 
\end{array}
\eqno(A.3)
$$ 

 $$
 \begin{array}{l} 
 {\ds \int _{0}^{\infty }J_{0} (2\sqrt{a t} ) \, k_{2 n} (t)\, dt=} 
 \\ \phantom{\rule{1pt}{15pt}}
  {\ds \frac{(-1)^{n\, -\, 1} }{2} \left[(n-1) ki_{2 n\, -\, 2} (a)-2 n ki_{2 n\, } (a)+(n+1) ki_{2 n\, +\, 2} (a)\right]} 
  \end{array}   
  \eqno (A.4)
$$

 $$
 \int _{0}^{\infty }J_{0} (2\sqrt{a t} ) \, k_{2 n} (t)\, \frac{dt}{t} =(-1)^{n} ki_{2 n} (a)
         \eqno (A.5)
         $$

 $$
 \int _{0}^{\infty }e^{-\, t} \, J_{1}  (2^{3/2} \sqrt{x t} )\, k_{2 n} (t)\, \frac{dt}{t} =\frac{(-1)^{n\, -\, 1} x^{n\, -\, 1/2} \, e^{-\, x} }{\sqrt{2}  n!} 
     \eqno(A.6) 
     $$

$$
\int _{0}^{\infty }e^{-\, a\, t} \, t^{n\, +\, 1/2} J_{1}  (2\sqrt{x t} )\, dt=\frac{(-1)^{n}  \Gamma (n+2)\, e^{-\, x/2 a} }{a^{n\, +\, 1} \, \sqrt{x} } \, k_{2 n\, +\, 2} \left(\frac{x}{2 a} \right)\, ;\quad a>0
\eqno (A.7)
$$

 $$
 \int _{0}^{x}\frac{J_{n} (t)-k_{2 n} (t)}{t}  \, dt=Ji_{n} (x)-ki_{2 n} (x)+\frac{(-1)^{n} }{n}         \eqno (A.8)
$$

 $$
 \int _{0}^{\infty }t^{n/2\, -\, 1} e^{-\, t} \, J_{2-n}  (4\sqrt{x t} )\, k_{2 n} (t)\, \frac{dt}{t} =\frac{x^{n/2\, -\, 1} \, e^{-\, x} }{2} k_{2 n} (x)
 \eqno (A.9)
$$

 $$
 \begin{array}{l} 
 {\ds \int _{0}^{\infty }\frac{e^{-\, b t^{2} } J_{\lambda } (a  \sqrt{t^{2} +x^{2} }  )
 \; J_{\nu } (a  \sqrt{t^{2} +x^{2} } ) }{t\, (t^{2} +x^{2} )^{(\lambda \, +\, \nu )/2} }
  \, k_{2 n\, +\, 2} (b t^{2} ) \, dt=\frac{(-1)^{n} J_{\lambda } (a x)\, J_{\nu } (a x)}{(2 n+2)\, x^{\lambda \, +\, \nu } } } 
 \\  \phantom{\rule{1pt}{15pt}}
 {\ds Re(\lambda +\nu )>-{3}/{2} }
  \end{array}
  \eqno(A.10)
$$

 $$
 \begin{array}{l}
  {\ds \int _{0}^{\pi }U_{2 n}  (\sqrt{x} \cos \theta )\, (\sin \theta )^{2} \, d\theta =
  \frac{\pi  (2 n)!\, e^{x/2} }{2 x n!} \, k_{2 n\, +\, 2} (\frac{x}{2} )}
   \\ \phantom{\rule{1pt}{15pt}}
    {\ds U_{n} (x)=(-1)^{n} \, e^{x^{2} } \, \frac{d^{n} }{d x^{n} }
     \left\{e^{-\, x^{2} } \right\}} 
    \end{array}
    \eqno (A.11)
    $$

 $$
 \int _{0}^{x}\sin (x-t) \, ki_{2} (t)\, dt=\cos x-\sin x-e^{-\, x} 
 \eqno(A.12)
 $$

 $$
 \int _{0}^{x}\cos (x-t) \, ki_{2} (t)\, dt=\cos x-\sin x+e^{-\, x} 
 \eqno (A.13)
 $$

 $$
 \int _{0}^{x}\sinh (x-t) \, ki_{2} (t)\, dt=e^{-\, x} (1+x)-\cosh x
 \eqno (A.14)
 $$

 $$
 \int _{0}^{x}\cosh (x-t) \, ki_{2} (t)\, dt=-x e^{-\, x} -\sinh x
 \eqno (A.15)
 $$

 $$
 \int _{0}^{x}e^{x\, -\, t}  \, ki_{2} (t)\, dt=-2 \sinh x
 \eqno (A.16)
$$

 $$
 \int _{0}^{x}(x-t)\, e^{x\, -\, t}  \, ki_{4} (t)\, dt=\sinh x-x \cosh x
 \eqno  (A.17)
 $$

 $$
 \int _{0}^{\infty }e^{\, -\, a t}  \, ki_{0} (b t)\, dt=\frac{1}{a} \, \ln \left(\frac{b}{a+b} \right)
 \, ;\quad a,b>0
 \eqno       (A.18)
 $$

$$
 \int _{0}^{\infty }\frac{ki_{2} (a t)-ki_{2} (b t)}{t}  \, dt=2\, \ln \left(\frac{a}{b} \right)
 \, ;\quad a,b>0
 \eqno (A.19)
$$

$$
\int _{0}^{\infty }J_{0} (2\sqrt{a t} ) \, ki_{2 n} (t)\, dt=(-1)^{n} \frac{k_{2 n} (a)}{a} 
 \eqno (A.20)
$$
\noindent

\section*{Appendix B: Trigonometric Integrals and Differential 
\\  Equations associated with the Kummer Confluent \\ Hypergeometric Functions}

\noindent In a paper devoted to the  note of Kummer, where he introduced into mathematics the confluent hypergeometric function defined by the following polynomial series

$$
 \begin{array}{l} 
 {\ds _{1} F_{1} (a;b;x)
 =M(a,b,x)
 =1+\frac{a }{b }  \frac{x}{1!} +\frac{a (a+1)}{b (b+1)}  \frac{x^{2} }{2!} +\frac{a (a+1) (a+2)}{b (b+1) (b+2)}  \frac{x^{3} }{3!} +\ldots }
  \\ \phantom{\rule{1pt}{15pt}}
   {\ds Re(b)>Re (a)>0}. 
   \end{array}
   \eqno(B.1)
   $$
the Italian mathematician  Giulio Giuliani [43] in 1888 considered the trigonometric integral (the original notation is replaced here by that used in this survey) 

$$
 {\ds I(x)=\int _{0}^{\pi /2}(\cos \theta )^{\alpha  - 1} \cos (\frac{x}{2} \tan \theta +n \theta )\, d\theta  =\frac{\pi }{2} \, k_{-n,\alpha -1} (\frac{x}{2} )\,,
 \quad {\alpha >1} }.
 \eqno(B.2)
$$
We note that  this integral is one of particular solutions of the following differential equation
$$
 {\ds 4  x \frac{d^{2} I(x)}{d x^{2} } -4 (\alpha -1) \frac{d I(x)}{d x}  -(x+2 n) I(x)=0}.
 \eqno(B.3)
$$
 Besides, Giuliani introduced two integrals coming from (B.2)
$$
 \begin{array}{l} 
 {\ds U_{n} (\alpha ,x)=\int _{0}^{\pi /2}(\cos \theta )^{\alpha - 1} \cos (\frac{x}{2} 
 \tan \theta )\cos (n \theta )\, d\theta  },
  \\ \phantom{\rule{1pt}{25pt}}
   {\ds V_{n} (\alpha ,x)=
  \int _{0}^{\pi /2}(\cos \theta )^{\alpha \, -\, 1} \sin (\frac{x}{2} \tan \theta )\sin (n \theta )
   \, d\theta  }, 
   \end{array}
   \eqno((B.4)
$$
when
$$
{\ds \int _{0}^{\pi /2}(\cos \theta )^{\alpha \, -\, 1} \cos (\frac{x}{2} \tan \theta +n \theta )
\, d\theta  
=U_{n} (\alpha ,x)-V_{n} (\alpha ,x)}.
\eqno (B.5)
$$
He showed that these integrals are solutions of the set of differential equations of the first order 
$$
 \begin{array}{l}
  {\ds 2 (\alpha -1)  \frac{d U_{n} (\alpha ,x)}{d x} +\frac{x}{2}   U_{n} (\alpha -2,x)  -  nV_{n} (\alpha ,x)  =0}, 
  \\  \phantom{\rule{1pt}{20pt}}
  {\ds 2 (\alpha -1)  \frac{d V_{n} (\alpha ,x)}{d x} +\frac{x}{2}   V_{n} (\alpha -2,x)  -  nU_{n} (\alpha ,x)  =0}, 
  \end{array}      
  \eqno(B.6)
  $$
and of the second order 
\vs
$$
 \begin{array}{l} 
 {\ds 2  x \frac{d^{2} U_{n} (\alpha ,x)}{d x^{2} } -2 (\alpha -1) 
 \frac{d U_{n} (\alpha ,x)}{d x}  -\frac{x}{2} \, U_{n} (\alpha ,x)+n  V_{n} (\alpha ,x)=0},
 \\  \phantom{\rule{1pt}{20pt}}
 {\ds 2  x \frac{d^{2} V_{n} (\alpha ,x)}{d x^{2} } -2 (\alpha -1) \frac{d V_{n} (\alpha ,x)}
 {d x}  -\frac{x}{2} \, V_{n} (\alpha ,x)+n  U_{n} (\alpha ,x)=0}. 
 \end{array}
   \eqno(B.7)
$$
\vs
 From (B.6) and (B.7) it is possible to obtain a differential equation of the fourth order
$$
\begin{array}{l} 
{\ds 4  x^{2}   \frac{d^{4} U_{n} (\alpha ,x)}{d x^{4} } -8 (\alpha -2) x  \frac{d ^{3} U_{n} (\alpha ,x)}{d x^{3} }  -}
 \\  \phantom{\rule{1pt}{20pt}} 
 {\ds 2  \left[x^{2} -2 (\alpha -1) (\alpha -2)\right]\, \frac{d ^{2} U_{n} (\alpha ,x)}{d x^{2} } +2 x (\alpha -2)  \frac{d U_{n} (\alpha ,x)}{d x} -}
  \\  \phantom{\rule{1pt}{20pt}} 
  {\ds \left(\frac{x^{2} }{4} +n^{2} +1-\alpha \right)\, U_{n} (\alpha ,x)=0}, 
  \\  \phantom{\rule{1pt}{20pt}} 
  {\ds V_{n} (\alpha ,x)=\frac{1}{n} \left(-\, 2 x  \frac{d ^{2} U_{n} (\alpha ,x)}{d x^{2} } +2 (\alpha -1)  \frac{d U_{n} (\alpha ,x)}{d x} +\frac{x}{2} \, U_{n} (\alpha ,x)\right)}. \end{array}
  \eqno(B.8)
$$
\vs
In terms of the Kummer confluent hypergeometric functions Giuliani was able to obtain that 
     $$
 \begin{array}{l} 
 {\ds \int _{0}^{\pi /2}(\cos \theta )^{\alpha \, -\, 1} \cos (\frac{x}{2} \tan \theta +n \theta )
 \, d\theta  =U_{n} (\alpha ,x)-V_{n} (\alpha ,x)=} 
 \\  \phantom{\rule{1pt}{25pt}}
  {\ds \left[\frac{\pi\, \Gamma (\alpha -1)\,  e^{-\, x/2} }{2^{\alpha } \, \Gamma \left(\frac{\alpha -n+1}{2} \right)\, \Gamma \left(\frac{\alpha +n+1}{2} \right)} \, _{1} F_{1} (\frac{\alpha -n+1}{2} ;1-\alpha ;x)-\right. }
   \\   \phantom{\rule{1pt}{25pt}}
   {\ds \left. \frac{\, \pi^2 \cos \left(\frac{\alpha -n}{2} \right)\, x^{\alpha } \, e^{-\, x/2} }{2^{\alpha } \, \sin (\pi  \alpha )\, \Gamma (\alpha )} \, _{1} F_{1} (\frac{\alpha +n+1}{2} ;\alpha +1;x)\right]},
       \end{array}
    \eqno (B.9)
    $$
and
$$
 \begin{array}{l} 
 {\ds U_{n} (\alpha ,x)+V_{n} (\alpha ,x)=\left[\frac{\pi \, \Gamma (\alpha -1) \, e^{-\, x/2} }{2^{\alpha } \, \Gamma \left(\frac{\alpha +n+1}{2} \right)\,
  \Gamma \left(\frac{\alpha -n+1}{2} \right)} \, _{1} F_{1} (\frac{1-\alpha -n}{2} ;1-\alpha ;x)\right. }, 
  \\  \phantom{\rule{1pt}{25pt}}
   {\ds \left. -\, \frac{\, \pi^2  \cos \left(\frac{\alpha +n}{2} \right)\, x^{\alpha } \, e^{-\, x/2} }{2^{\alpha } \, \sin (\pi  \alpha )\, \Gamma (\alpha )} \, _{1} F_{1} (\frac{\alpha -n+1}{2} ;\alpha +1;x)\right]}. 
  \end{array}
   \eqno(B.10)
$$
\noindent These expressions can be presented in terms of the generalized Bateman functions defined in \eqref{GrindEQ__75_}
$$
 \begin{array}{l} 
 {\ds k_{-\, \nu ,\alpha ,0} (x)=} 
 \\  \phantom{\rule{1pt}{25pt}}
 {\ds \left[\frac{    \Gamma (\alpha )\,     e^{-\, x} }{2^{\alpha } \, \Gamma \left(\frac{\alpha -\nu }{2} +1\right)\,
  \Gamma \left(\frac{\alpha +\nu }{2} \right)} \, 
   _{1}F_{1} (\frac{\alpha -\nu }{2} +1;-\alpha ;2    x)-\right. } 
   \\  \phantom{\rule{1pt}{25pt}}
    {\ds \left. \frac{\pi \,     \cos \left(\frac{\alpha -\nu +1}{2}
     \right)\, x^{\alpha \, +\, 1} \, e^{-\, x} }{2^{\alpha } \, \sin [\pi (    \alpha +1)]\, \Gamma (\alpha +1)} \, _{1}F_{1} (\frac{\alpha +\nu }{2} +1;\alpha +2;2    x)\right]}, 
    \end{array}
    \eqno (B.11)
    $$
\noindent and
$$
 \begin{array}{l} 
 {\ds k_{\nu ,\alpha ,0} (x)=\left[\frac{\, \Gamma (\alpha )        e^{-\, x} }{2^{\alpha } \, \Gamma \left(\frac{\alpha +\nu }{2} +1\right)\, \Gamma \left(\frac{\alpha -\nu }{2} +1\right)} \, _{1}F_{1} (\frac{-\alpha -\nu }{2} ;-\alpha ;2    x)\right. }
  \\\phantom{\rule{1pt}{25pt}}
   {\ds \left. -\, \frac{\pi \,     \cos \left(\frac{\alpha +\nu +1}{2} \right)\,
    x^{\alpha \, +\, 1} \, e^{-\, x} }{2^{\alpha } \, \sin [\pi     (\alpha +1)]\, \Gamma (\alpha +1)} \, _{1}F_{1} (\frac{\alpha -\nu }{2} +1;\alpha +2;2    x)\right]}.
    \end{array}
     \eqno(B.12)
$$

 As shown by Giuliani, by changing the integration variable, the finite trigonometric integrals can be presented as the infinite integrals, for example
$$
 {\ds \int _{0}^{\pi /2}(\cos \theta )^{\alpha } \cos (\frac{x}{2} \tan \theta )\, d\theta  =
 \int _{0}^{\infty }\frac{\cos \left(\frac{x t}{2} \right)}{(1+t^{2} )^{\alpha /2 + 1} }
 \, dt\,.}
   \eqno(B.13)
$$
 Considering the case $\alpha=1$ in (B-2), Bateman  [44] in 1931 noted the link that exists between the investigated by Giuliani integral and the $k$-Bateman function with negative order. 
 He also found that the solution of the following third order differential equation
$$
 {\ds x \frac{d^{3} I(x)}{d x^{3} } -  (\alpha -1) \frac{d^{2}  I(x)}{d x^{2} }  
 -(x+ n) \frac{d I(x)}{d x} -\beta  I(x)=0,}   
 \eqno(B.14)
$$
 is given by the following trigonometric integral 
$$
 {\ds I(x)=\, \int _{0}^{\pi /2}(\cos \theta )^{\alpha } (\sin \theta )^{\beta  - 1}
  \cos (x\tan \theta +n \theta )\, d\theta 
   =\frac{\pi }{2} \, k_{n,\alpha ,\beta  - 1} (x)\,.}
   \eqno(B.15)
$$
Besides, Bateman showed that for $x>0$:
$$
 \begin{array}{l} 
 {\ds \int _{0}^{\pi /2}(\cos \theta )^{m} \cos [x\tan \theta +(m+2n) \theta ]\, d\theta  =\frac{e^{x} \sin (\pi  n)}{2^{k\, +\, 1} } \, \int _{0}^{1}t^{k}  (1-t)^{n\, -\, 1} 
 \, e^{-\, 2 x/t} \, dt\,,} 
 \\   \phantom{\rule{1pt}{15pt}}
 {\ds \int _{0}^{\pi /2}\cos [x\tan \theta +(m+2n) \theta ]\, d\theta  =\frac{e^{x} \sin (\pi  n)}{2} \, \int _{0}^{1}(1-t)^{n\, -\, 1} \, e^{-\, 2 x/t}  \, dt}
 ={\ds \frac{\pi }{2} \, k_{-\, 2 n} (x)\,.} 
 \end{array}
 \eqno(B.16)
$$

As can be observed, the included material from  the 1888 paper by Giuliani and 
 from the 1931 paper by Bateman 
is important from the historical and mathematical points of view.

\section*{Appendix C: Integral Representations  of Special Functions used in this Survey}

\noindent 
\textbf {Hypergeometric Function}      
$$
{\ds      _2 F_ 1      (a,b;c;x) =
\frac{\Gamma (c)}{\Gamma (a)\,\Gamma (b)}
                   \, \int _0^1    \frac  {t^ {b- 1}\, (1-t)^{c -b- 1}}{ (1-xt)^a} \,dt \,
                   \quad  Re(c)>Re(a)>0\,.}                            
                     \eqno(C.1)
$$

\noindent 
\textbf {Kummer Confluent Hypergeometric Function}

$$
{\ds      _1 F_ 1      (a,b;x) = M(a,b,x) =
\frac{\Gamma (b)}{\Gamma (a)\,\Gamma (b-a)}
                   \, \int _0^1    t^ {a- 1}\, e^{xt}\, (1-t)^{b -a- 1} \,dt \,
                   \quad  Re(b)>Re(a)>0\,.}                            
                     \eqno(C.2)
$$

\noindent \textbf{Tricomi Confluent Hypergeometric Function}     

$$
    {\ds U(a,b,x)=\frac{1}{\Gamma (a)}\, \int _0^\infty  t^ {a - 1} \,e^{-xt}\,(1+t)^{b-a-1}\,dt\,,                 
     \quad  Re(b)>Re(a)>0\,.}                            
                     \eqno(C.3)
     $$

\noindent 
\textbf{Whittaker Functions}     
$$
\begin{array}{l}
{\ds    M_{\kappa ,\mu}(x)
=\frac{\Gamma (1+2\mu )\, x^{\mu  + 1/2}\,     e^{- x/2}}
{\Gamma (\mu +\kappa +1/2)\, \Gamma (\mu -\kappa +1/2)}
     \,\int_0^1  t^{\mu - \kappa  - 1/2} \,  e^{ x t}\, (1-t)^{\mu  + \kappa  - 1/2}\, dt\,,}
     \\  \phantom{\rule{1pt}{15pt}}
      {\ds M_{\kappa ,\mu}(x)
      =x^{\mu  + 1/2}      e^{ - x/2}\,   M(\mu -\kappa + 1/ 2,1+2 \mu ,x)}
       \\  \phantom{\rule{1pt}{15pt}}
         {\ds      Re(\mu \pm \kappa + 1/2)>0\,.}
            \end {array}    
               \eqno(C. 4)
               $$

\noindent 

$$
\begin{array}{l}
{\ds    W_{\kappa ,\mu}(x)
=\frac{x^{\mu  + 1/2}\,     e^{- x/2}}
{ \Gamma (\mu -\kappa +1/2)}
     \,\int_0^\infty t^{\mu - \kappa  - 1/2} \,  e^{- x t}\, (1+t)^{\mu  + \kappa  - 1/2}\, dt\,,}
     \\  \phantom{\rule{1pt}{15pt}}
      {\ds W_{\kappa ,\mu}(x)
      =x^{\mu  + 1/2}      e^{ - x/2}\,   U(\mu -\kappa + 1/ 2,1+2 \mu ,x)}
       \\  \phantom{\rule{1pt}{15pt}}
         {\ds      Re(\mu - \kappa + 1 2)>0\,.}
            \end {array}    
               \eqno(C.5)
               $$

\noindent \textbf{Bessel Functions}     

\noindent 
$$
{\ds J_\nu  (x)=
\frac{1}{\pi}   \, \int _0^ \pi  \cos (x \sin \theta -\nu \theta)\, d \theta -
\frac{\sin (\pi  \nu )}{\pi}\, \int_0 ^ \infty  e ^{ -x  \sinh t -\nu t}\,  dt.}
\eqno(C.6)
$$
\noindent 

$$
{\ds Y_\nu  (x)=
\frac{1}{\pi}   \, \int _0^\pi  \sin (x \sin \theta -\nu \theta)\, d \theta -
\frac{\sin (\pi  \nu )}{\pi}\, \int_0 ^\infty  e ^{ -x  \sinh t -\nu t}
\, \left[  e^{\nu t} + e^{-\nu t} \cos (\pi \nu)\right]\,  dt.}
\eqno(C.7)
$$

$$
{\ds I_\nu  (x)=
\frac{1}{\pi}   \, \int _0^\pi 
e^{x \cos \theta}  \cos (\nu \theta)\, d \theta -
\frac{\sin (\pi  \nu )}{\pi}\, \int_0 ^\infty  e^{ -x  \cosh t -\nu t}\,  dt.}
\eqno(C.8)
$$

$$
{\ds K_\nu  (x)=
\frac{\Gamma(\nu+1/2)\, (2x)^\nu}{\sqrt \pi} 
  \, \int _0^\infty 
\frac{\cos(xt)}{(1+t^2)^{\nu+1/2}}\, d t
= \int_0^\infty  e^{ -x  \cosh t}\, \cosh (\nu t)\,  dt.}
\eqno(C.9)
$$

\noindent \textbf{Struve Functions}     

$$
{\ds H_\nu (x)
=\frac{2 (x/2)^\nu} {\Gamma (\nu +1/2)   \sqrt  \pi}
    \,  \int_0^1     (1-t^ 2)^{\nu  - 1/2}\,  \sin (x t )\, dt\,, \quad Re(\nu)> -1/2.}
    \eqno(C.10)
$$

$$
{\ds L_\nu (x)
=\frac{2 (x/2)^\nu} {\Gamma (\nu +1/2)   \sqrt  \pi}
    \,  \int_0^{\pi/2}      (\sin  t )^{2\nu}\, \sinh(x \cos t)\, dt\,, \quad Re(\nu)> -1/2.}
    \eqno(C.11)
$$

\noindent \textbf{Lommel Functions}     

$$
{\ds S_{\mu, \nu} (x)
x^ \mu      \, 
\int_0^\infty      e^{- x  t}   \, _2 F_1 
     \left(  \frac{1-\mu +|nu} {2}, \frac{1-\mu -\nu}{ 2} ; \frac{1}{ 2}; -t ^ 2  \right)
     \, dt\,, \quad    Re(x)>0\,.}
             \eqno (C. 12)
   $$


\end{document}